\documentclass[11pt]{article}
\usepackage[pctex32]{graphics}
\usepackage{amssymb}
\usepackage{latexsym}
\usepackage{epsfig}
\usepackage{multirow}
\usepackage{pstcol}
\usepackage{float}
\usepackage{epstopdf}
\usepackage{mathrsfs}
\usepackage{graphicx}
\usepackage{amsmath}
\usepackage{amssymb}
\usepackage{relsize}
\usepackage{multirow}

\usepackage{graphicx}
\usepackage{multicol}

\usepackage{times}
\usepackage{epstopdf}

\newtheorem{definition}{Definition}[section]
\newtheorem{theorem}[definition]{Theorem}

\definecolor{Red}{rgb}{1,0.,0.}
\usepackage{fancyhdr}
\usepackage{setspace}

\newcommand{\R}{{\mathbb R}}

\newcommand{\mA}{{\mathsf A}}
\newcommand{\mB}{{\mathsf B}}
\newcommand{\mC}{{\mathsf C}}
\newcommand{\mU}{{\mathsf U}}
\newcommand{\mV}{{\mathsf V}}
\newcommand{\mX}{{\mathsf X}}
\newcommand{\mL}{{\mathsf L}}
\newcommand{\mD}{{\mathsf D}}

\newcommand{\mP}{{\mathsf P}}

\newcommand{\mK}{{\mathsf K}}
\newcommand{\mY}{{\mathsf Y}}

\newcommand{\mSigma}{{\mathsf \Sigma}}
\newcommand{\mT}{{\mathsf T}}
\newcommand{\mI}{{\mathsf I}}
\newcommand{\mzero}{{\mathsf 0}}

\begin{document}

\title{Bayes meets Krylov: preconditioning CGLS\\ for underdetermined systems}
\author{D Calvetti$^1$ \and F Pitolli$^2$ \and E Somersalo$^1$ \and B Vantaggi$^2$}
\date{$^1$Case Western Reserve University \\ Department of Mathematics, Applied Mathematics and Statistics \\ 10900 Euclid Avenue, Cleveland, OH 44106 \\
$^2$ University of Rome ``La Sapienza'' \\ Department of Basic and Applied Science for Engineering \\ via Scarpa 16, 00161 Rome, Italy }
\maketitle

\begin{abstract}
The solution of linear inverse problems when the unknown parameters outnumber data requires addressing the problem of a nontrivial null space. After restating the problem within the Bayesian framework, a priori information about the unknown can be utilized for determining the null space contribution to the solution. More specifically, if the solution of the associated linear system is computed by the Conjugate Gradient for Least Squares (CGLS) method, the additional information can be encoded in the form of a right preconditioner. In this paper we study how the right preconditioned changes the Krylov subspaces where the CGLS iterates live, and draw a tighter connection between Bayesian inference and Krylov subspace methods. The advantages of a Krylov-meet-Bayes approach to the solution of underdetermined linear inverse problems is illustrated with computed examples.

\end{abstract}


\section{Introduction}
Severely underdetermined inverse problems in which the dimensions of the unknown greatly outnumber the data are commonly encountered in applications, and the challenge is to introduce in a numerically efficient way additional and complementary information to overcome the problem of scarcity of the data. The Bayesian statistical framework provides a systematic way to augment the observation model with prior information \cite{CSbook,KSbook,Stuart,Tarantola}; however, the implementation of effective algorithms for computing informative pointwise estimators pose a challenge.
In this paper we are interested in solving severely underdetermined linear inverse problems with the Conjugate Gradient for Least Squares (CGLS) method equipped with a suitable stopping rule. In particular, we want to analyze how the Krylov subspace method is affected when we address the non-uniqueness of the solution within the Bayesian statistical framework.  This work is, in particular, motivated by the interest in large-scale underdetermined linear inverse problems that arise in important biomedical applications, such as electrical impedance tomography  (EIT) \cite{Cheney} and magnetoencephalography (MEG) \cite{Baillet}, in which the CGLS method combined with Bayesian models has proven to be particularly efficient, see, e.g. \cite{DebraThesis,McGivney2,Homa}.

To introduce the setting, assume that we want to estimate a vector $x\in\R^n$ from the observation of the vector $b\in\R^m$ satisfying
\begin{equation}\label{b=fx}
 b = F(x) + \varepsilon,\quad F:\R^n\to\R^m,
\end{equation}
where $m<n$, and $\varepsilon$ represents additive noise. In the Bayesian setting all unknowns are modeled as random variables \cite{CSbook} described by their probability distributions that we assume here to be expressed in terms of probability density functions (pdf). Denote by $\pi_{\rm noise}$ the pdf of the noise and let the prior density $\pi_{\rm prior}$  encode the information about $x$ prior to considering the observations. If $x$ and $\varepsilon$ are statistically independent, the posterior density of $x$ conditioned on the observed value of $b$ is, according to Bayes' formula, given by
\begin{equation}\label{bayes}
 \pi(x\mid b) \propto \pi_{\rm prior}(x)\pi_{\rm noise}(b - F(x)),\quad b =b_{\rm observed},
\end{equation}
where ``$\propto$'' stands for ``proportional up to a normalizing constant''. Here we concentrate on observation models that are liner, that is, $F(x) = \mA x$, where $\mA\in\R^{m\times n}$ is a given matrix. The approach that we analyze carries over to non-linear problems through a sequence of local linearizations, as has been demonstrated in several articles: see, e.g., \cite{ArridgeLSQR,McGivney,DebraThesis,McGivney2,Pursiainen,Pursiainen2}. For the sake of clarity, we limit our discussion to the case of Gaussian distributions. The results of our analysis can be extend naturally to hierarchical Bayesian models and furthermore to non-Gaussian problems,  through distributions belonging to exponential family; however, the extension will be addressed in a separate contribution.

An important application that motivates this work and gives an idea of the dimensions of the problems that we have in mind is the MEG inverse problem, in which the aim is to estimate electric activity of the brain by measuring the weak magnetic fields outside the skull. This is a classical inverse source problem for Maxwell's equations and can be described in terms of a linear mapping: $x\in\R^n$ represents the discretized sources inside the head, $b\in\R^m$ are the magnetometer recordings outside the head, and $\mA\in\R^{m\times n}$ is the lead field matrix. In a typical MEG problem, $m$ represents the number of channels in the device, and is of the order $m \sim 100-300$, while the number of unknowns is typically orders of magnitude larger,  from tens of thousands upwards. In particular, since $m\ll n$, a significant amount of additional information is needed to find a reasonable solution.

In this paper we assume that the system (\ref{b=fx}) of linear equation that arise from the discretization of the underlying linear inverse problem is solved using a Krylov subspace iterative method, and in order that such a solution is representative in view of the Bayesian model (\ref{bayes}),
additional prior information is embedded in the algorithm via a right preconditioning matrix. If the linear system is not square, as is the case for the problems of interest to us, the Krylov subspace iterative solver of choice is CGLS, originally proposed in \cite{HS}. The analysis of the effects of a statistically inspired right preconditioner on the Krylov subspaces is carried out in detail here only for the CGLS method, but analogous arguments and techniques can be used to extend the analysis to other Krylov subspace iterative methods.

Originally preconditioners for linear systems were introduced to increase the convergence rate of iterative methods, a feature that made them not well suited for the solution of linear discrete inverse problems, because the speed-up may also increase the rate at which the amplified noise components contaminate the solution.  In fact, a fast converging CGLS iteration often returns a solution that satisfies the data adequately without capturing important features of the solution.
To overcome this problem, special preconditioners for ill-posed problems were proposed in, e.g., \cite{HNP,Enrich2003}. The idea in this class of regularizing preconditioners was to accelerate the convergence rate of the portion of the solution in the subspace of the signal, while leaving the noise subspace unpreconditioned, a task that could be achieved with a preconditioning matrix approximating the inverse of $\mA$ in the signal subspace and acting as an identity on the noise subspace. Constructing such matrix would require a priori knowledge of the signal and noise subspaces. In the case of underdetermined ill-posed systems, preconditioning can be used to enrich the computed solution, albeit at the cost of slowing down the rate of convergence of the iterative method.

In general, the solution of a linear system can be decomposed as the sum of a {\em visible} component which contributes to the data, and an {\em invisible} component, that belongs to the null space of the matrix. Effectively, the CGLS method equipped with a suitably chosen right preconditioner can leverage rich prior information about the solution to extract a significant component of the solution from the effective null space of the matrix $\mA$. The coupling between the visible and invisible part of the solution happens through non-orthogonal subspace decompositions based on the prior information implicitly applied by the iteration method. One of the aims of this work is to understand and quantify the changes to the fundamental subspaces of the linear system induced by the right prior conditioners, information that can be subsequently used in the design and selection of statistically inspired right preconditioners.

The paper is organized as follows. In Section~\ref{sec:linear}, we review some basic facts on linear inverse problems and regularization, including the standard and preconditioned CGLS as a regularization method. Using subspace decompsitions, we then analyze how the priorconditioning alters the fundamental subspaces and the induced effect on the solutions. We also discuss the convergence rate of the priorconditioned  system through the Lanczos process. Finally in Section~\ref{sec:examples}, we elucidate the ideas with two computed examples.

\section{Linear inverse problems: Tikhonov vs. Krylov-subspace regularization}\label{sec:linear}

Consider the discretized linear version of (\ref{b=fx}),
\begin{equation}
\label{linsys_noise}
 b = \mA x + \varepsilon,
\end{equation}
where $\mA\in\R^{m\times n}$, $m<n$ and $\varepsilon$ represent the noise due, e.g.,  measurement errors and model uncertainties. Adhering to the paradigm of Bayesian inference, all unknowns are modeled as random variables and are described in terms of their probably density functions.  The belief about the unknown $x$ prior to taking the data into consideration is expressed by its prior density, that here we assume to be a Gaussian density with mean zero and symmetric positive definite covariance matrix  $\mC\in\R^{n\times n}$. Furthermore, we assume that the noise vector $\varepsilon$  is a zero mean white Gaussian random variable.  More general Gaussian noise can be reduced to this case multiplying both sides of the linear system by a suitable invertible matrix, a process known in signal processing as whitening. It follows from Bayes' formula (\ref{posterior}) that the posterior density of $x$, which is the solution of the inverse problem in the Bayesian setting, is given by
\begin{equation}\label{posterior}
 \pi(x\mid b) \propto {\rm exp}\left( -\frac 12 \|\mA x - b\|^2  - \frac 12 x^\mT \mC^{-1} x\right).
\end{equation}
Consider a symmetric decomposition of the inverse of the covariance matrix, or the precision matrix, of the form
\begin{equation}\label{factoring}
 \mC^{-1} = \mB^\mT \mB;
\end{equation}
whose existence is guaranteed by the positive definiteness of $ \mC$; here $\mB$ can be, e.g., a Cholesky factor, or the square root of the precision matrix.
By using the factorization (\ref{factoring}), the negative logarithm of the posterior density (\ref{posterior}), also known as Gibbs energy, can be written as
\[
 G(x) = \|\mA x - b\|^2 + \| \mB x\|^2 = \left\|\left[\begin{array}{c} \mA \\ \mB\end{array}\right]x - \left[\begin{array}{c} b \\ 0\end{array}\right]\right\|^2,
\]
from which it follows that the maximizer of  (\ref{posterior}), known as the maximum a posteriori (MAP) estimate of $x$, can be found by solving the system
\[
\left[\begin{array}{c} \mA \\ \mB\end{array}\right]x = \left[\begin{array}{c} b \\ 0\end{array}\right]
\]
in the least squares sense. The MAP estimate is also the solution of the linear system
\[
 \big(\mA^\mT \mA + \mB^\mT\mB\big)x =  \mA^\mT b,
\]
which are the normal equations for the penalized least squares problem
\begin{equation}\label{Tikhonov}
x = {\rm argmin} \{ \| b - \mA x\|^2 + \lambda \|\mB x \| ^2 \}
\end{equation}
associated with Tikhonov regularization with linear regularization function $\mB$ and regularization parameter $\lambda$ equal to unity \cite{CSbook}. The standard version of Tikhonov regularization assumes that $\mB = \mI_n$, the unit matrix. The solution of (\ref{Tikhonov}) can be computed by first transforming the problem into standard form, an operation that requires the generalized singular value decomposition (GSVD) of the matrix pair $(\mA,\mB)$: see, e.g., \cite{hansenGSVD}, as we will review below. Strategies to reduce the high computational costs associated with the computation of the GSVD in the context of Tikhonov regularization have been proposed in \cite{Dykes}. Several numerical methods for the computation of Tikhonov regularized solution for large problems can be found  in the literature; see, e.g., \cite{CMRS,Tikhonov2003,Tiklarge}. The determination of a suitable value for the regularization parameter for Tikhonov in standard form has been studied extensively in the literature, where it has been related to the amount of error in the right hand side. In the case where Tikhonov regularization operator is different from the identity it is not obvious how to determine the value of the regularization parameter without transforming the problem to standard form.

When the unknown $x$ has a large number of components, a computationally attractive alternative to Tikhonov regularization with regularization operator $\mB$ is to use a Krylov subspace iterative methods, equipped with a suitable stopping rule, to solve approximately the linear system
\begin{equation}\label{linsys}
 \mA x = b,
\end{equation}
with the matrix $\mB$ as a right preconditioner \cite{Priorconditioners2,Priorconditioners} as explained below. In the case where the matrix $\mA$ is non square, a natural choice is  to use the  CGLS method. The CGLS method determines a sequence of approximate solutions of the linear equations associated with the linear system (\ref{linsys}) without explicitly forming the matrix $\mA^\mT\mA$, but instead multiplying vectors with $\mA$ and $\mA^\mT$ separately.

Iterative solution methods for large scale non square linear systems have been studied extensively, and different implementations of the idea behind the CGLS method which take into account the characteristics of the problem have been proposed in the literature. For a discussion of the different implementations and guidelines their relative advantages see, e.g.,  \cite{Bjorck,PaigeSaunders,ChoiSaunders,Saad,vdV2003}.

The stopping rule is an important component of iterative linear systems solvers. In general, the iteration stops when the norm of residual error is sufficiently small. It has been shown that when the linear systems is the discretization of a linear ill-posed problems, the CGLS method with a suitably modified stopping rule is a regularization method. In the next subsection we recall a few results about Krylov subspace regularization for underdetermined linear systems, and we compare the regularized solutions computed by  CGLS in standard form and CGLS with a right preconditioner.

\subsection{Standard and priorconditioned CGLS}

The regularizing properties of the CGLS method equipped with a suitable stopping rule are well known: see, e.g., \cite{hanke,hankehansen,HNP,hansenbook,hansenbook2}.  Most of the studies of the properties of CGLS method as a regularization scheme have been carried out for overdetermined linear systems, where the number of equations exceeds the degrees of freedom of the problem. While a lot of the results carry over to the case of underdetermined problems, we will see below how the presence of a null space, which in some cases may be very large, changes the picture.

The CGLS method for solving the linear system (\ref{linsys}) starting with the initial approximate solution $x_0=0$ computes a sequence of approximate solutions $x_1, \ldots, x_k$ until a termination criterion is satisfied. If the matrix $\mA$ is ill-conditioned and the right side noisy, as is typically the case when (\ref{linsys}) is the discretization of a linear inverse problem, the iterations are stopped as soon as the norm of the discrepancy, $d_k = b - \mA x_k$ falls below a threshold value corresponding to the noise level in the data. The  CGLS iterates are computed by projecting (\ref{linsys}) onto a nested family of Krylov subspaces. More precisely,  the $j$th iterate satisfies
\[
x_j={\rm arg min}\big\{ \|\mA x - b\|\mid x\in  {\mathbf{\mathcal{K}}}_j(\mA^\mT b, \mA^\mT \mA)\big\}
\]
where
 \begin{equation}\label{krylov}
{\mathcal K}_j( \mA^\mT b, \mA^\mT \mA)= {\rm span} \{ \mA^\mT b, (\mA^\mT \mA)\mA^\mT b, \ldots, (\mA^\mT \mA)^{j-1}\mA^\mT b \}
\end{equation}
is the $j$th Krylov subspace associated with the method. From the assumption that the noise $\varepsilon\in\R^m$ is additive, zero mean white Gaussian it follows that
\[
 {\mathsf E}\big\{\|\varepsilon\|^2\big\} = m.
\]
If the rank of the Krylov subspaces continues to increase, the iteration is stopped as soon as
\[
 \|\mA x_k - b\|^2 <\tau m,
\]
where $\tau>0$ is a safeguard factor, e.g., $\tau = 1.2$, and $x_k$ is the regularized solution determined by the method. Typically the stopping index $k$ is much smaller than $m$. This last observations implies that when using the CGLS method for the solution of discrete inverse problems, it is important that the salient features of the solution vector are included  in the first few iterates already, because it may require only a few iteration steps to satisfy the stopping criterion.

In general, the presence of a non-trivial null space for the matrix $\mA$ raises the need to address how to deal with the degrees of freedom in the solution that cannot be determined by the data. The standard CGLS method sets the null-space components of the computed solutions to zero, seeking to find the solution of minimum Euclidean norm. This corresponds to look for the maximum likelihood estimator.

Within the Bayesian framework, setting the null space contribution to the solution to zero is only justified if there are reasons to believe that the solution is orthogonal to the null space. If that is not the case, the way to proceed is to encode a priori information in the prior density and have it determine  the contribution of the null space to the computed solution. The encoding of the information carried by the prior into the CGLS iterations via a right preconditioned was formalized in \cite{Priorconditioners2} and in \cite{Priorconditioners}, where the term {\em prior conditioner} was first proposed, although the use of smoothing right preconditioners for linear discrete ill-posed problems has also been advocated earlier in, e.g., \cite{invertibleSmoothing,hanke,hankehansen}.

Without loss of generality, we can assume a priori that $x\sim{\mathbf {\mathcal N}}(0,\mC)$. If $\mB$  is the matrix defined in (\ref{factoring}), the change of variables
\[
 w = \mB x \sim {\mathbf {\mathcal N}}(0,\mI_n)
\]
transforms $x$ into a zero mean white Gaussian random variable, hence we can write the linear system (\ref{linsys}) in terms of the new variable as
\begin{equation}\label{linsys white}
 \mA\mB^{-1} w = b,
\end{equation}
where
\begin{equation}\label{w to x}
x = \mB^{-1}w.
\end{equation}
The CGLS algorithm applied to the system (\ref{linsys white}) with starting vector $w_0=0$ computes a sequence $w_0, w_1, \ldots$ of approximate solutions for (\ref{linsys white}) and a corresponding sequence of vectors in the original coordinate system through (\ref{w to x}).  It turns out that this simple change of variables determines a sequence of approximate solution which may be very different from that computed by the CGLS method applied directly to (\ref{linsys}). To investigate the effect of the priorconditioning on the Krylov subspaces where the approximate solutions are computed in relation to the null space of the coefficient matrix, we need some linear algebra results.

Consider the canonical orthogonal decomposition of the space of the solution into the null space of $\mA$ and its orthogonal complement, which is the range of $\mA^\mT$,
\begin{equation}\label{null space}
 \R^n = {\mathscr N}(\mA)\oplus {\mathscr R}(\mA^\mT).
\end{equation}
It follows from (\ref{null space}) that the Krylov subspace (\ref{krylov}) is orthogonal to ${\mathscr N}(\mA)$, hence $x_k \perp {\mathscr N}(\mA)$.  In other words, all the CGLS iterates  are perpendicular to the null space of $\mA$.

Introducing the notation $\widetilde\mA = \mA\mB^{-1}$ to simplify the notation,  we observe that
the  $j$th iterate of the whitened problem (\ref{linsys white}) solves the minimization problem
\[
w_j ={\rm arg min}\big\{ \|\widetilde \mA w - b\|\mid w\in  {\mathcal{K}}_j(\widetilde\mA^\mT b, \widetilde\mA^\mT \widetilde\mA)\big\},
\]
and the corresponding $j$th {\em Priorconditioned CGLS (PCGLS) solution}  $\widetilde x_j = \mB^{-1}w_j$ satisfies
\begin{equation}\label{p_krylov}
 \widetilde x_j\in{\rm span}\big\{\mB^{-1}(\widetilde A^\mT \widetilde A)^\ell \widetilde A^\mT b\mid  0\leq \ell\leq j-1\big\}.
\end{equation}
It follows from
\[
 \mB^{-1}\widetilde A^\mT = \mB^{-1}\mB^{-\mT}\mA^\mT = \mC\mA^\mT,
\]
that
\begin{equation}\label{p_krylov2}
\mB^{-1}\big(\widetilde A^\mT\widetilde\mA\big)^\ell\widetilde \mA^\mT = \big(\mC\mA^\mT\mA)^\ell\mC\mA^\mT,\quad 0\leq \ell\leq j-1.
\end{equation}
Therefore, in view of (\ref{p_krylov})--(\ref{p_krylov2}),
\begin{equation}\label{xj}
 \widetilde x_j \in \mC \big({\mathscr N}(\mA)^\perp\big),
\end{equation}
and $\widetilde x_j$ is not necessarily orthogonal to the null space of $\mA$. This observation is discussed further in the light of subspace factorizations in the next section.

\subsection{Subspace factorizations}

To characterize the subspaces where the iterates of the PCGLS method are computed, we need some results about a generalization of the singular value decomposition.

The  generalized singular value decomposition (GSVD), first introduced in \cite{VanLoan} and then described  in \cite{GolubVanLoan}, is a standard tool to analyze regularized solutions produced by Tikhonov regularization with a regularization functional different from the identity \cite{hansenbook,hansenbook2,hansenGSVD,Elden}.
In our analysis we resort to formulation proposed in \cite{PaigeSaundersGSVD}.
\begin{theorem}
Given a matrix pair $(\mA,\mB)$\, with $\mA \in \R^{m\times n}$, $\mB \in \R^{n\times n}$, $m<n$, it is possible to find a factorization of the form
\[
 \mA = \mU\left[ \begin{array}{cc} \mzero_{m,n-m} & \mSigma_\mA\end{array}\right] \mX^{-1},\quad \mB = \mV \left[\begin{array}{cc} \mI_{n-m} & \\ & \mSigma_\mB\end{array}\right]\mX^{-1},
\]
called the generalized singular value decomposition, where $\mU\in\R^{m\times m}$ and $\mV\in\R^{n\times n}$ are orthogonal matrices, $\mX\in\R^{n\times n}$ is an invertible matrix, and $\mSigma_\mA\in\R^{m\times m}$ and $\mSigma_\mB\in\R^{m\times m}$ are diagonal matrices.
The diagonal entries $s^{(\mA)}_1, \ldots, s^{(\mA)}_m$ and $s^{(\mB)}_1, \ldots, s^{(\mB)}_m$ of the matrices $\mSigma_\mA$ and $\mSigma_\mB$ are real, nonnegative and satisfy
\begin{eqnarray}
s^{(\mA)}_1& \leq & s^{(\mA)}_2 \leq \ldots \leq  s^{(\mA)}_m \nonumber \\
s^{(\mB)}_1 & \geq & s^{(\mB)}_2 \geq \ldots \geq  s^{(\mB)}_m \nonumber \\
(s^{(\mA)}_j)^2 &+& (s^{(\mB)}_j)^2 = 1, \qquad 1 \leq j \leq m. \label{sumsquare}
\end{eqnarray}
 It follows from  the condition (\ref{sumsquare}) that $0 < s^{(\mA)}_j \leq 1$ and $ 0 < s^{(\mB)}_j \leq 1$.
The ratios $s^{(\mA)}_j/s^{(\mB)}_j$ for $1 \leq j \leq m$ are called the generalized singular values of the matrix pair $(\mA,\mB)$.
\end{theorem}

If $\mA$ has full row rank, as we assume for simplicity here, the diagonal entries of $\mSigma_\mA$ are positive.

Given a symmetric positive definite matrix $\mC$, the $\mC$-inner product, or prior energy inner product, of two vectors $x$ and $y$ is defined as
\[
 \langle x,y\rangle_\mC = x^\mT \mC^{-1} y.
\]
Two vectors $x$ and $y$ are  $\mC$-orthogonal, denoted by $ x\perp_\mC y $ if and only if $
 \langle x,y\rangle_\mC = 0$. When $x,y$ are random variables $\mC$-orthogonality corresponds to independence of the transformed variables $\mB x$ and $\mB y$ (where the matrix $\mB$ is defined in (\ref{factoring})).

The following theorem sheds some light on the structure of the null space of $\mA$ in the GSVD basis. Given a matrix $\mY$, we use the notation ${\rm span}(\mY)$ to denote the subspace spanned by the columns of the matrix $\mY$.

\begin{theorem} If we partition the matrix $\mX\in\R^{n\times n}$ in Theorem 2.1 as
\[
 \mX = \left[\begin{array}{cc} \mX' & \mX''\end{array}\right], \quad \mX' \in\R^{n\times (n-m)},\; \mX''\in\R^{n\times m},
\]
it follows that
\[
 {\mathscr N}(\mA) = {\rm span}\big\{ \mX'\big\},
\]
and we can express $\R^n$ as a $\mC$-orthogonal direct sum,
\[
 \R^n = {\rm span}\big\{\mX'\big\} \oplus_\mC {\rm span}\big\{\mX''\big\} = {\mathscr N}(\mA)\oplus_\mC {\rm span}\big\{\mX''\big\}.
\]
\end{theorem}

{\em Proof.} Let $z\in{\mathscr N}(\mA)$ and define
\[
 \zeta = \mX^{-1} z =\left[\begin{array}{c} \zeta ' \\ \zeta ''\end{array}\right], \quad \zeta'\in\R^m,\; \zeta'' \in\R^{n-m}.
\]
From the observation that
\[
 \mA z  = \mU \left[ \begin{array}{cc} \mzero_{m,n-m} & \mSigma_\mA\end{array}\right]\zeta = \mU \mSigma_\mA\zeta '' =0,
\]
the orthogonality of $\mU$, and the positivity of the diagonal elements in $\mSigma_\mA$  it follows that $\zeta'' = 0$. This implies that
\[
 z = \mX\zeta = \mX' \zeta'\in {\rm span}\big(\mX'\big),
\]
that is,
\[
 {\mathscr N}(\mA)\subset {\rm span}\big(\mX'\big).
\]
It is easy to see, by following the argument backwards, that the inclusion holds in the other direction also.

To show the $\mC$-orthogonal decomposition of the space $\R^n$ in terms of the basis $\mX$, observe that
\[
 \mC^{-1}  = \mB^\mT\mB = \mX^{-\mT} \left[\begin{array}{cc} \mI_{n-m} & \\ & \mSigma_\mB^2\end{array}\right] \mX^{-1},
\]
or, equivalently,
\begin{equation}\label{C matrix}
 \mC  = \big(\mB^\mT\mB\big)^{-1}  = \mX \left[\begin{array}{cc} \mI_{n-m} & \\ & \mSigma_\mB^{-2}\end{array}\right] \mX^\mT,
\end{equation}
which is tantamount  to saying  that $\mC$ is diagonal in the basis $\mX$.

To show that
\[
 {\rm span}\big\{\mX'\big\} \perp_\mC {\rm span}\big\{\mX''\big\},
\]
let $x\in {\rm span}\big\{\mX'\big\}$ and $y\in {\rm span}\big\{\mX''\big\}$. Then
\[
 x = \mX'\alpha = \mX\left[\begin{array}{c}\alpha \\ 0\end{array}\right],\quad y = \mX''\beta = \mX \left[\begin{array}{cc} 0\\ \beta\end{array}\right]
\]
 therefore
\[
 \langle x,y\rangle_\mC = \left[\begin{array}{cc} \alpha^\mT & 0\end{array}\right]\mX^\mT \mX^{-\mT}\left[\begin{array}{cc} \mI_{n-m} & \\ & \mSigma_\mB^2\end{array}\right]
\mX^{-1}\mX \left[\begin{array}{cc} 0\\ \beta\end{array}\right] =  \left[\begin{array}{cc} \alpha^\mT & 0\end{array}\right]\left[\begin{array}{cc} 0\\
\mSigma_\mB^2\beta\end{array}\right] =0,
\]
as claimed. Hence, any vector $x\in\R^n$ can be expressed in the $\mX$-basis as
\[
 x = \mX\left[\begin{array}{c} \alpha \\ \beta\end{array}\right] = \mX'\alpha + \mX''\beta
\]
in a unique way, thus completing the proof.  $\Box$

If follows immediately from the theorem above that
\[
 {\mathscr N}(\mA)^\perp = {\mathscr R}(\mA^\mT),\quad  {\mathscr N}(\mA)^{\perp_\mC} = {\rm span}\big\{\mX''\big\}.
\]
This factorization allows us to interpret the effect of priorconditioning: According to (\ref{xj}), we have
\[
 \widetilde x_j \in\mC\big( {\mathscr R}(\mA^\mT)\big).
\]
In particular, if ${\mathscr R}(\mA^\mT)$ is an invariant subspace of the covariance matrix $\mC$, then the iterates $\widetilde x_j$ are orthogonal to the null space of $\mA$, and the PCGLS is not capable of informing the null space. In particular, in statistical terms, if the orthogonal projections of $x$  onto the subspaces ${\mathscr N}(\mA)$ and ${\mathscr R}(\mA^\mT)$ are uncorrelated, then these spaces are $\mC$-orthogonal, which implies that $\mC({\mathscr R}(\mA^\mT)) \subset {\mathscr R}(\mA^\mT)$. Conversely, if the projections are correlated, the component of the PCGLS iterates may have a significant component in the null space of $\mA$. This lack of $\mC$-orthogonality is the key to include a priori information about the contribution from the null space, hence invisible to the likelihood, into the computed solutions.

A good measure for the level of correlation introduced by the priorconditioner may therefore be the distance from  $\mC$-orthogonality of the null space of $\mA$ and the range of its transpose. This distance can be quantified, e.g., in terms of the $\mC$-inner product as
\[
 {\mathscr I}(\mA,\mC) = \inf \left\{\frac{|\langle x,y\rangle_\mC|}{\|x\|_\mC\|y\|_\mC} \,\mid\, x\in {\mathscr N}(\mA),\; y\in{\mathscr R}(\mA^\mT)\right\}.
\]

Observe that the subspace decompositions above show that the $\mC$-orthogonal projections of $x$ on ${\mathscr N}(\mA)$ and ${\rm span}\big\{\mX''\big\}$ are always uncorrelated;  see also  \cite{hansenOblique} for discussion of oblique projections in the context of standard transformations.

\subsection{The Lanczos process}

The convergence rate of the CGLS algorithm is related to the spectral properties of the matrix, and consequently, it is of interest to understand how preconditioning changes the spectral properties and the order in which eigendirections are included in the computed solution by the CGLS method. We consider the Lanczos process implicitly defined by the Krylov subspace iterations.
It is known \cite{MMQ} that the first $k$ residual vectors computed by CGLS and normalized to have unit length form an orthonormal basis for the Krylov subspace $\mathcal{K}_k ( \mA^\mT b,\mA^\mT \mA )$. Denote them by $v_0, v_1, \ldots, v_{k-1} $ and collected them into the matrix $\mV_k$. If
\[ x_j=x_{j-1}+ \alpha_{j-1}p_{j-1}
\]
is the CGLS updating formula which computes the $j$th iterate $x_j$ from the $(j-1)$st by adding a correction along the search direction $p_{j-1}$,  and
\[ p_j=r_j+ \beta_{j-1}p_{j-1} \]
is the formula to update the search direction with $r_j= b -\mA^T \mA x_j$ , it is straightforward to show that
\[
\mA^\mT \mA V_k = V_k \mT_k - \frac{\sqrt{\beta_{k-1}}}{\alpha_{k-1}} v_k e_k^\mT.
\]
In the last formula,  $e_k$ is the $k$th column of the $k \times k$ identity matrix and $\mT_k$ is the tridiagonal matrix
\[
\mT_k = \mL_k ^{-1}\Delta_k ^{-1}\mL_k,
\]
where $\Delta_k = {\rm diag} \{ \alpha_0, \ldots, \alpha_{k-1} \}$, $\mL_k = \Phi_k \mU_k \Phi_k^{-1}$,
 $\Phi_k= {\rm diag} \{ \|r_0,\|, \ldots, \|r_{k-1}\| \}$, and $\mU_k$ is the $k \times k$ upper bidiagonal matrix with ones on the main diagonal and $-\beta_0, \ldots, -\beta_{k-1}$ on the superdiagonal.

It follows from the orthogonality of the $v_\ell$ that
\[
\mV_k^\mT (\mA^\mT \mA) V_k = \mT_k,
\]
therefore $\mT_k$ is the projection of $\mA^\mT \mA$ onto the Krylov subspace ${\mathcal{K}}_k ( \mA^\mT b,\mA^\mT \mA )$. One can show that the $k$th CGLS iterate can be expressed as
\[
x_k =\mV_ky_k,
\]
where $y_k$ solves the $k \times k$ linear system
\[
\mT_k y = \|r_0\| e_1.
\]

Before proceeding with our analysis, we recall the following result about the rate of convergence of the conjugate gradient method.

\begin{theorem}
If the matrix $\mA$ is symmetric positive definite, and $\eta_k = x_*-x_k$ is the error in the $k$th iterate of the conjugate gradient method, then
\[
\| \eta_k\|_\mA \leq 2 \left( \frac{\sqrt{\kappa} -1}{\sqrt{\kappa}+1} \right)^k \|\eta_0\|_{\mA},
\]
where $\kappa = {\lambda_n}/{\lambda_1}$ is the condition number of $\mA$.
\end{theorem}

Replacing the matrix $\mA$ with $\mA^\mT\mA$, and noting that because of the termination criterion effectively $\lambda_n$ can be replaced by the smallest nonzero eigenvalue of $\mA^\mT \mA$, we can rewrite the theorem above as follows.

\begin{theorem}
If $x_k $ is the $k$th iterate computed with the CGLS method and the initial residual $\mA^\mT \mA $, it follows that
\[
\|r_k\| \leq 2 \left( \frac{\sqrt{\kappa} -1}{\sqrt{\kappa}+1} \right)^k \|r_0\|.
\]
where $\kappa$ is the condition number of $\mA^\mT\mA$.

\end{theorem}

Exploiting the connection between the Lanczos process, orthogonal polynomials and Gauss quadrature rules, we have the following characterization of the norm of the residual for the normal equations solved with the CGLS method in terms of the eigenvalues of the matrix $\mA^\mT\mA$ and of their approximation with the corresponding Ritz value, which are the eigenvalues of the tridiagonal matrix $\mT_k$. This is an adaptation of a result in  \cite{MMQ}  (page 203).

\begin{theorem}\label{discrepancy}
Let $\lambda_i\in\R$ denote the $i$th eigenvalue of $\mA^\mT\mA$, and $q_i\in\R^n$ be the corresponding eigenvector, $1\leq i\leq n$. Further, let $\theta_j^{(k)}$ be the $j$th eigenvalue of the tridiagonal matrix $\mT_k$, $1\leq j\leq k\leq r$, where $r$ is the rank of the matrix $\mA$.
For all $k$, $1 \leq k \leq r$,  there exists a number $\xi_k$,  $\lambda_1 \leq \xi_k \leq \lambda_r$ such that the norm of the $k$th residual satisfies
\[
\|r_k\|^2 = \frac{1}{\xi_k^{2k+1}} \sum_{i=1}^n \left[  \prod_{j=1}^k \left(\lambda_i - \theta_j^{(k)} \right)^2 \right] \left(r_0^{\mT}q_i \right)^2,
\]
\end{theorem}

In particular, the rate of convergence of CGLS depends on how accurately the eigenvalues $\lambda_i$ are approximated by the Ritz values $\theta_j^{(k)}$, and in particular those eigenvalues corresponding to eigendirections with a significant component of the initial residual $r_0$: As soon as a Ritz value has converged to $\lambda_i$, the contribution to the residual from the corresponding eigendirection vanishes. In turn, the convergence rate of the CGLS method depends on how accurately the eigenvalues of projected tridiagonal matrix approximate those of $\mA^\mT\mA$, see, e.g.,  \cite{vdSvdV}, hence it is reasonable to assume that the $k$th iterate is richer along those directions that approximate those eigenvectors of $\mA^\mT\mA$ where the components of $r_0$ are larger.

In general, if the vector $r_0$ is the linear combination of a few eigendirections of $\mA^\mT\mA$, in exact arithmetic it will require only a few iterations of the CGLS methods to meet the termination criterion, and the solution vectors will live in a low dimensional Krylov subspace methods. Conversely, the richer in eigendirections $r_0$ is, the more iteration will take the CGSL method to satisfy the stopping rule, and the larger the dimension of Krylov subspace of the solution will be.

It has been observed experimentally that the introduction of a prior conditioner in the CGLS method changes the spectral properties of the underlying normal equations in two ways. In fact, the projections of the initial residual vector along the eigenvectors of $\mA^\mT \mA$ associated with nonzero eigenvalues are more even distributed, and the range of the nonzero eigenvalues widens, sometimes significantly. In view of these observations and of Theorem~\ref{discrepancy}, we expect the convergence rate to slow down, because a larger number of eigenvalues of $\mA^\mT \mA$ must be approximated accurately by the Ritz values, in order to reduce the norm of the discrepancy below the assigned threshold. A consequence of this is that a larger number of eigendirections contribute to the computed solution,  with the result of providing a richer approximation of the underlying signal. We illustrate this in the next section with computed examples.

\section{Computed examples}\label{sec:examples}

In this section, we consider two linear inverse problems to elucidate the analysis of the PCGLS algorithm. The purpose of the first  simple one-dimensional test model is to follow step by step the changes in the Krylov subspaces and related quantities induced by the introduction of a prior conditioner, while the second example is a more realistic problem inspired by the medical imaging applications that motivated the present analysis.

\subsection{Example 1: One-dimensional deconvolution}

Let $f:[0,1]\to\R$ be a piecewise continuous function, $f(0) = 0$, and assume that the observations consist of few noisy convolutions by an Airy kernel function
\[
 a(t) = \left(\frac{J_1(\kappa t)}{\kappa t}\right)^2,
\]
where $\kappa>0$ is a parameter regulating the width of the kernel. To discretize the problem, we subdivide the interval $[0,1]$ by introducing $n$ equidistant points $s_1, \ldots, s_n$ and approximate the convolution integral by a finite sum of the form
\[
g(t) = \int_0^1 a(t - s) f(s) ds \approx \frac 1n \sum_{k=1}^{n} a(t-s_k) f(s_k),\quad 1\leq j\leq n,
\]
where $s_j = j/n$. We assume that the data consist of a $m$ discrete noisy measurements of $g$ at $t_1, \ldots, t_m$, that is
\[
 b_\ell = g(t_{\ell}) + \varepsilon_\ell,\quad 1\leq \ell \leq m,
\]
where $m \ll n$.  If we denote by $\mA$ the matrix whose entries are
\[
a_{j,k}= \frac 1n a(t_j - s_k), \qquad 1\leq j \leq m, \; 1 \leq k \leq n,
\]
and by $x$ the vector with components
\[
x_k= f(s_k), \qquad 1\leq k \leq n,
\]
we  can formulate the problem as the linear system (\ref{linsys_noise}).

In this computed example, we set $n = 150$ and $m=6$, thus the linear system is strongly underdetermined. We quantify the index of under determinacy by the ratio of unknowns to equation, which in this case is $n/m = 25$. The width parameter in the simulations is $\kappa =  0.02$. We solve this linear system iteratively starting with initial approximate solution $x_0=0$ in two different ways: first we use the plain CGLS method with a termination criterion based on the discrepancy principle, then we  solve it with the priorconditioned CGLS method. The prior conditioner corresponds to a second order Gaussian prior, and is defined through the factorization of its precision matrix,
\[
\mC^{-1} = \mL^\mT \mL,\quad  \mL = \beta \left[\begin{array}{rrrrr} \alpha  & & &  & \\ -1 & 2 & -1 & & \\ &  & \ddots & & \\ & & -1& 2 & -1 \\ & & & & \alpha \end{array}\right],
\]
where the value of the scalar $\alpha>0$ is selected to guarantee uniform variance for all pixels (see \cite{CSbook}).

\begin{figure}
\centerline{\includegraphics[width=12cm,height =6cm]{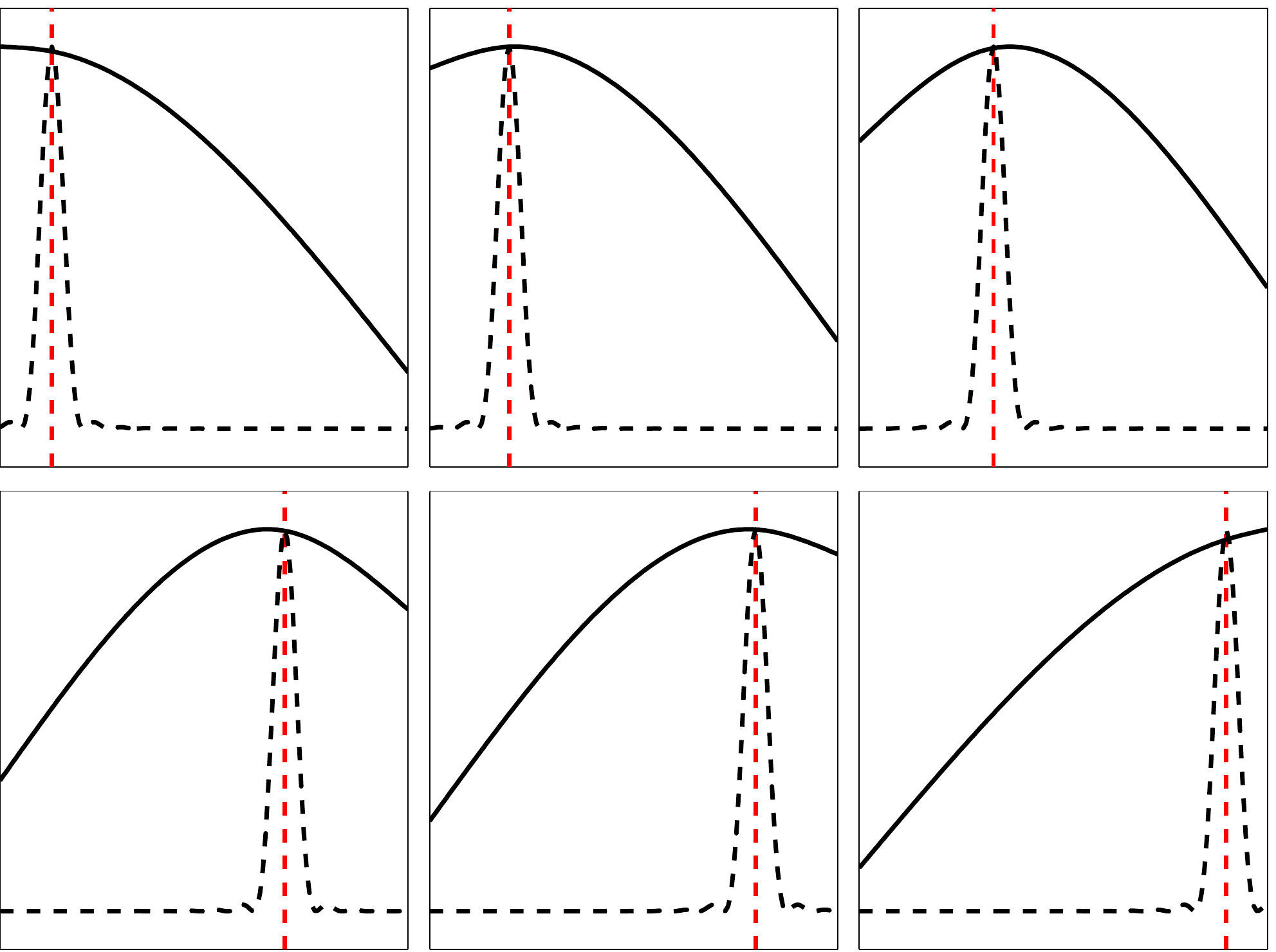}}
\caption{\label{fig:basis} The six basis vectors that span ${\mathscr R}(\mA^\mT)$ (dashed curve), and the vectors that span $\mC \big({\mathscr R}(\mA^\mT)\big)$, respectively (solid curve). The values of $t_j$ corresponding to the data are indicated by the vertical red lines.}
\end{figure}

The six basis functions spanning the subspace ${\mathscr N}(\mA)^\perp={\mathscr R}(\mA^\mT)$ used by the plain CGLS to compute the iterates, and the six vectors spanning $\mC\big({\mathscr N}(\mA)^\perp\big)=\mC\big({\mathscr R}(\mA^\mT)\big)$ are show in Figure~\ref{fig:basis}.  As expected, the basis functions reflect well the expected behavior of the underlying function as postulated by the chosen prior.

\begin{figure}
\centerline{\includegraphics[width=7cm]{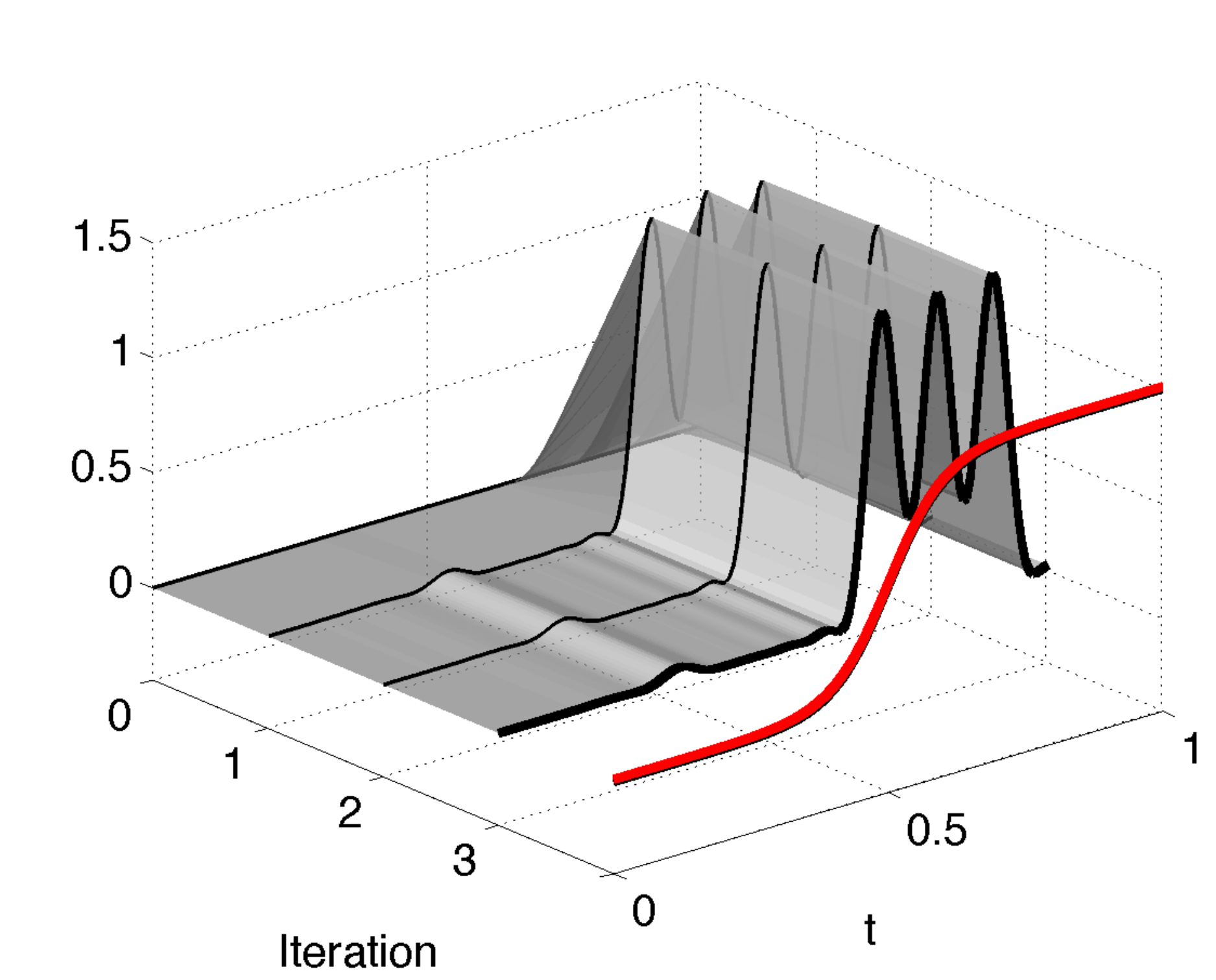}
\includegraphics[width=7cm]{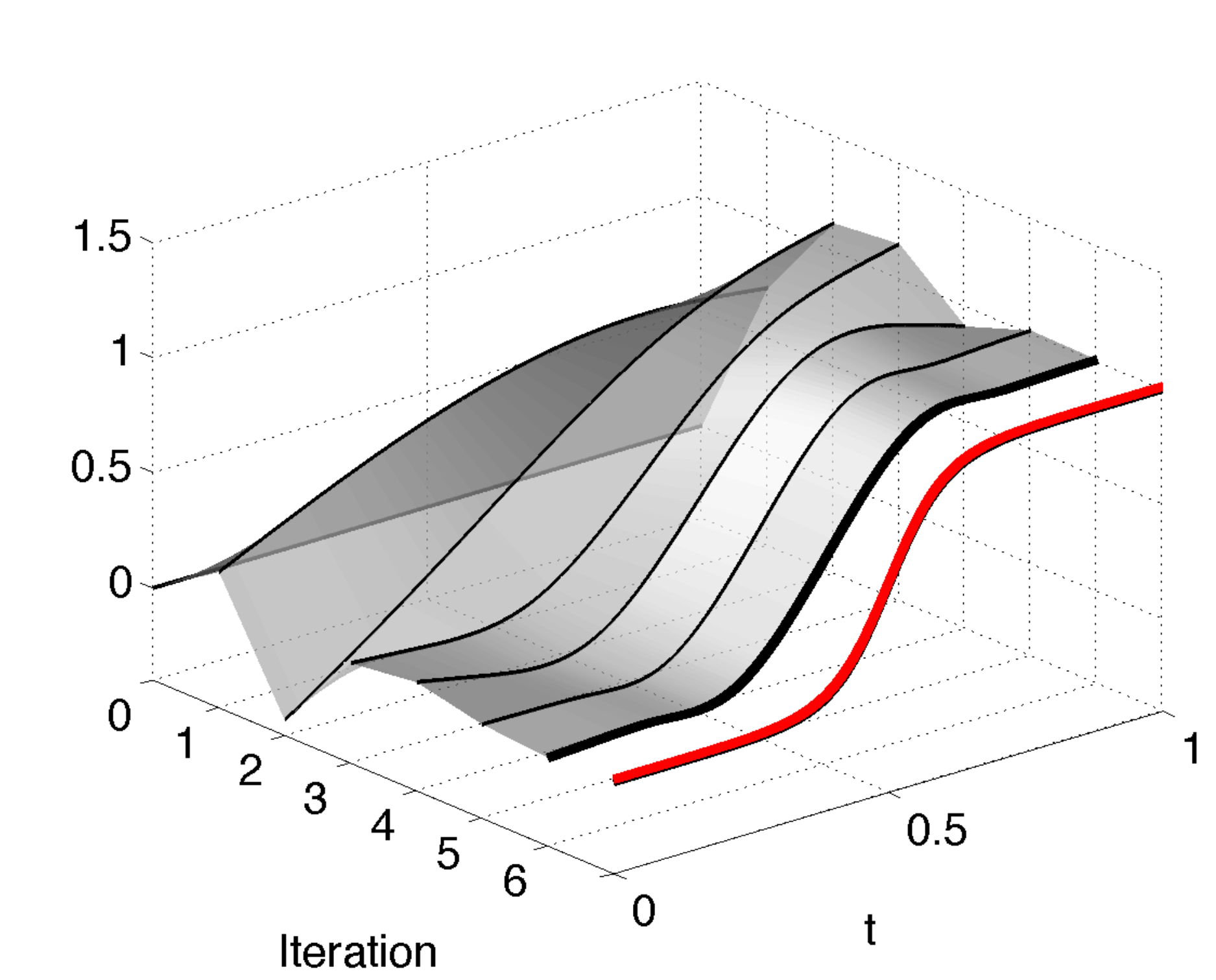}
}
\caption{\label{fig:deconv results}
The progress of the iterations based on plain CGLS (left) and the priorconditioned CGLS (right). The final result when iterations are stopped at discrepancy is the boundary shape of the surface. The true profile that was used to generate the data is shown in red.}
\end{figure}

To demonstrate the effect of the priorconditioning on the CGLS algorithm, we generate noisy data assuming that the underlying signal is a smooth sigmoid function and the additive noise in the data is scaled zero mean white Gaussian, $\varepsilon \sim {\mathbf {\mathcal N}}(0,\sigma^2\mI_m)$. The noise level is chosen very low, $\sigma = 5\times 10^{-5}$,  in order to allow the algorithm to compute more iterates; with higher noise level, the stopping criterion is satisfied after one or two iterations.
Figure~\ref{fig:deconv results} shows the sequence of CGLS and PCGLS iterates. In this example the standard CGLS converges very quickly, leading to a termination after three steps, while the PCGLS requires more iterations. Moreover, since the CGLS cannot inform the solution about the component of the solution in the null space of $\mA$, it returns a very poor reconstruction in the intervals between the observation points, as shown in the left plot of Figure~\ref{fig:deconv results}. The PCGLS algorithm, on the other hand, is able to extract meaningful information about the solution from  the null space of $\mA$ guided by the data and the prior, and returns a better solution.

\begin{figure}
\centerline{\includegraphics[width=6cm]{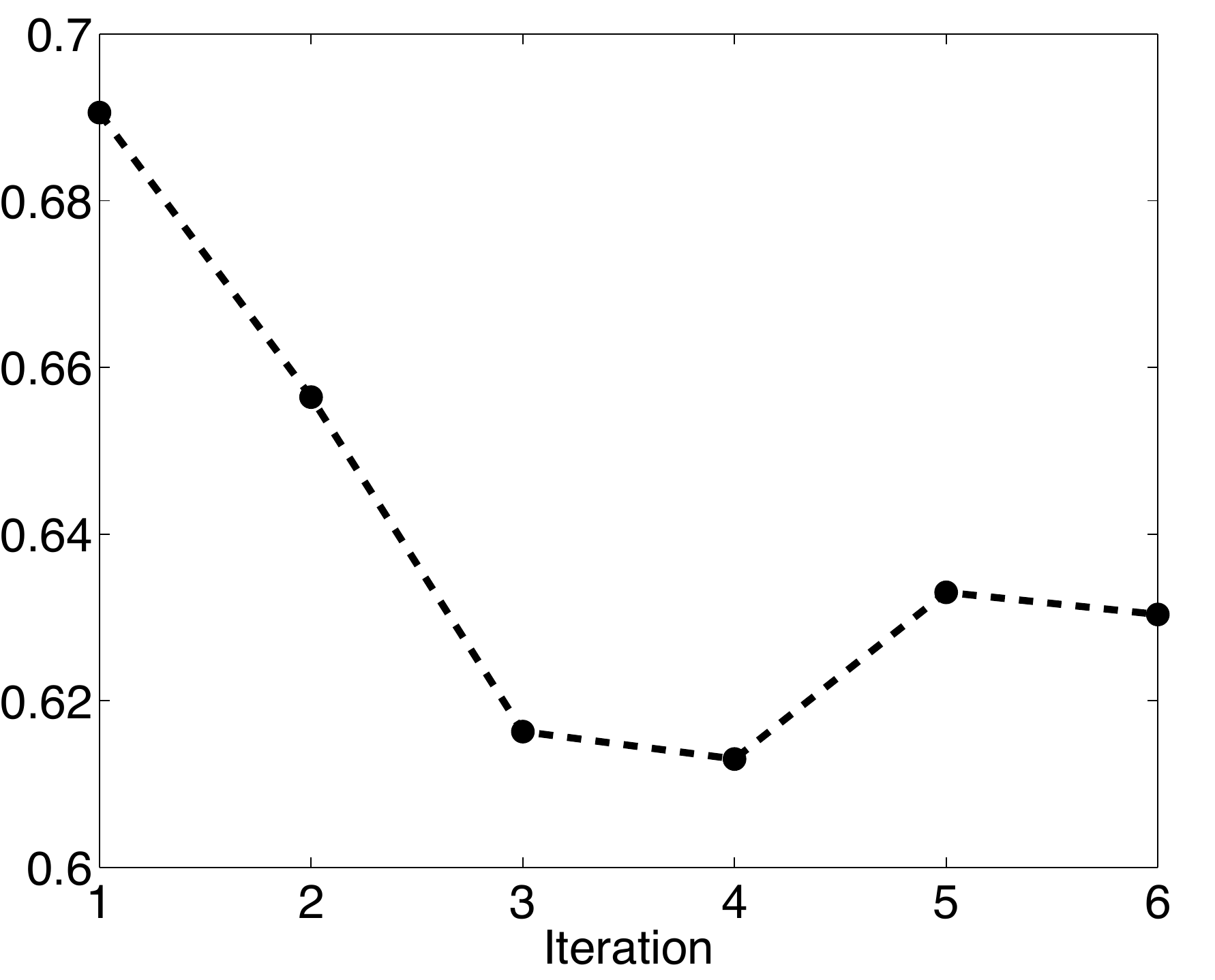}}
\caption{\label{fig:null space} The norm of the projection of the current iterate onto the null space divided by the norm of the current iterate, plotted as a function of the iteration number. The plot indicates that more than 60\% of the priorconditioned iterates is in the null space, while without the priorconditioning, the null space component vanishes.}
\end{figure}

An important role of the priorconditioner is to unlock important information about the underlying solution hidden in the null space by producing iterates not necessarily orthogonal to null space of $\mA$.
To quantify the distance from  $\mC$-orthogonality of the span of $\mX''$ and the null space of $\mA$ introduced by the priorconditioner at each iteration step,  we compute the relative size of the orthogonal component of the PCGLS iterates $\widetilde x_k$ on the null space of $\mA$.
\[
 \nu_k = \frac{\|\mP\widetilde x_k\|}{\|\widetilde x_k\|},\quad \mP:\R^n \mathop{\longrightarrow}^{\perp} {\mathscr N}(\mA).
\]
The plot of these components as a function of the iteration number, shown in Figure~\ref{fig:null space}, indicates that in this example more than 60\% of the priorconditioned solution is in the null space, while without the priorconditioning, the null space component vanishes.

\begin{figure}
\centerline{\includegraphics[width=6cm]{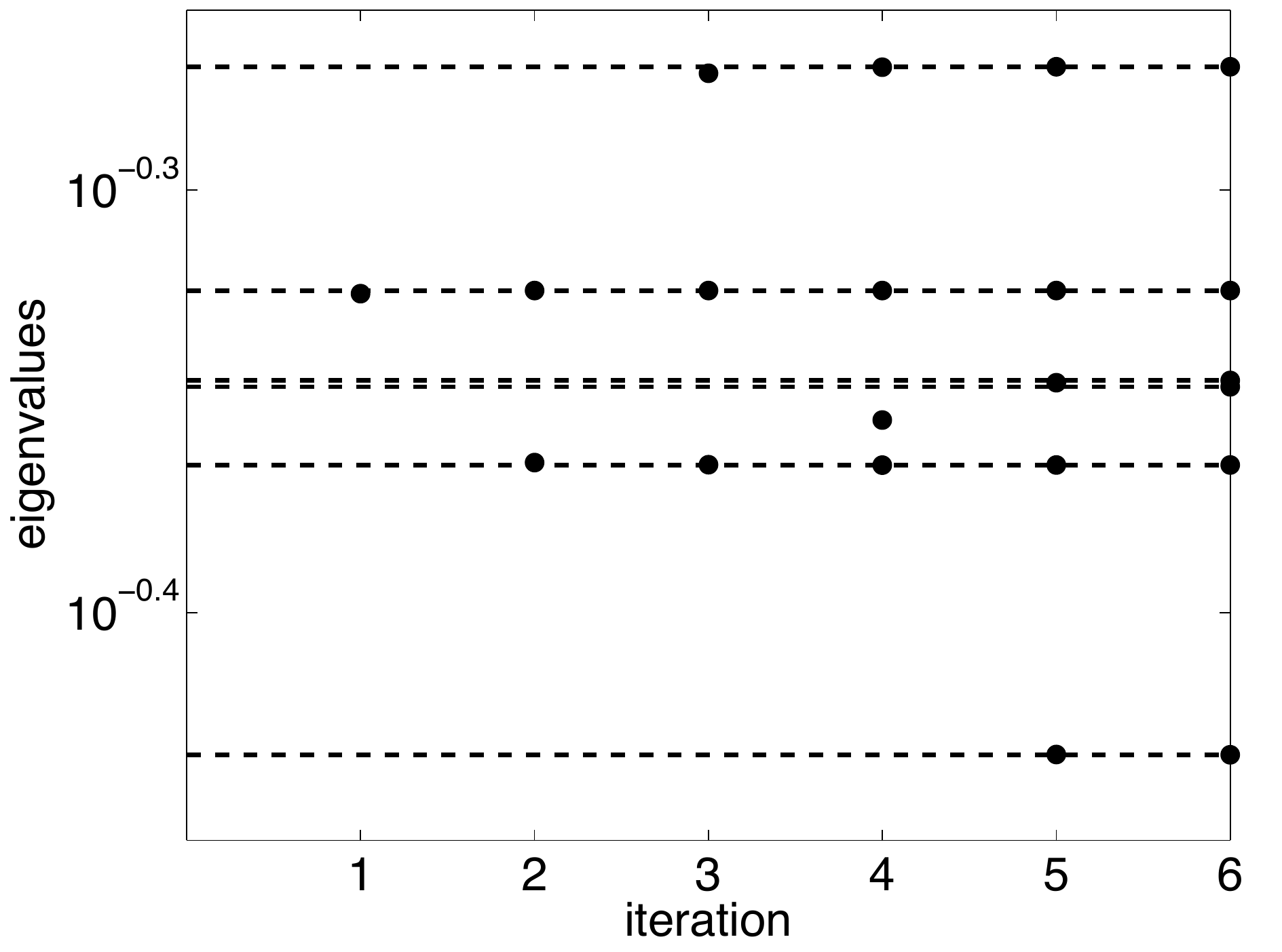}
\includegraphics[width=6cm]{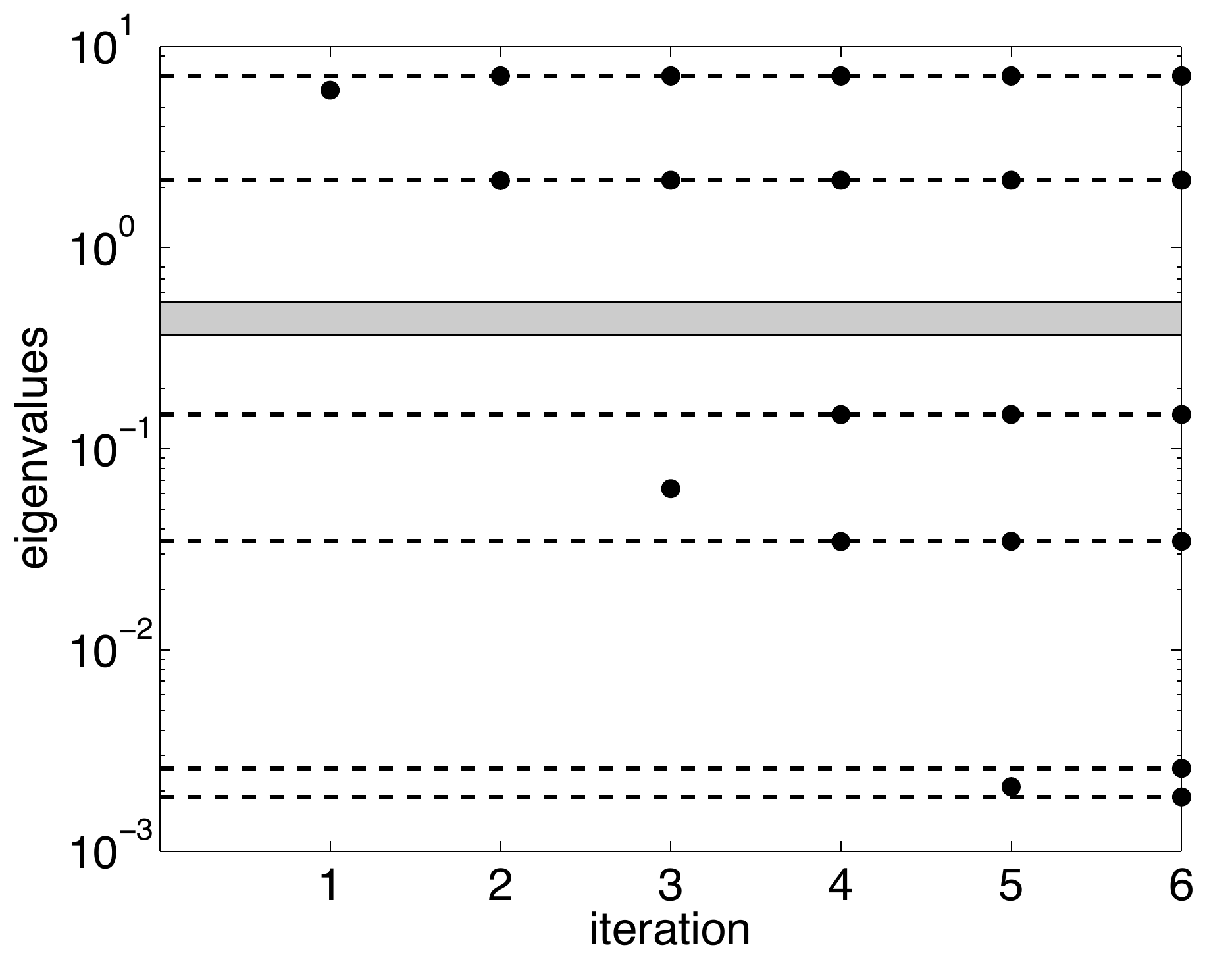}
}
\caption{\label{fig:spectrum deconv} Left: The eigenvalues of the tridiagonal matrix $\mT_k\in\R^{k\times k}$ as a function of the iteration number $k$. The non-zero eigenvalues of the matrix $\mA^\mT\mA$ are indicated by horizontal lines. Observe that two of the eigenvalues are nearly equal. Right: The corresponding eigenvalues using the priorconditioned scheme. To help comparing the two cases, the spectral interval of the matrix without priorconditioning is indicated as a shaded band in the plot.}
\end{figure}

To illustrate how this null-space unlocking  prior conditioner changes the eigenvalues of the least squares problem, we compute the spectra of the projected matrices
\[
 \mT_j = \mV_j^\mT (\mA^\mT \mA)\mV_j , \quad   \widetilde\mT_j = \widetilde\mV_j^\mT (\widetilde\mA^\mT \widetilde\mA)\widetilde\mV_j \quad 1\leq j \leq m,
\]
and the projections of the respective initial residuals onto the eigenvectors associated with the non-zero eigenvalues. Figure~\ref{fig:spectrum deconv} shows the non-zero eigenvalues and the progression of their approximation with the eigenvalues of the Lanczos tridiagonal matrices as $j$ increases.

\begin{table}[htdp]
\caption{Absolute values of he projections of the initial residual on the eigenvectors of $\mA^\mT \mA$ (left) and $\widetilde\mA^\mT\widetilde \mA$ (right). In agreement with Figure~\ref{fig:spectrum deconv}, the order in which the first four Ritz values (ORV) converge to the eigenvalues correspond the order of the largest projections.}
\begin{center}
\begin{tabular}{ cccc}
Plain CGLS & ORV & Preconditioned CGLS & ORV \\
\hline
0.0196 & (3) & 8.2586  & (1) \\
2.9481 & (1)  & 4.3337  & (2) \\
0.0329  & & 0.0759 &  (4) \\
0.0232 & (4) & 0.1214 &  (3) \\
0.4073 & (2) & 0.0031 & \\
0.0004 & & 0.0044 &
\end{tabular}
\end{center}
\label{default}
\end{table}

\begin{figure}
\centerline{\includegraphics[width=6cm]{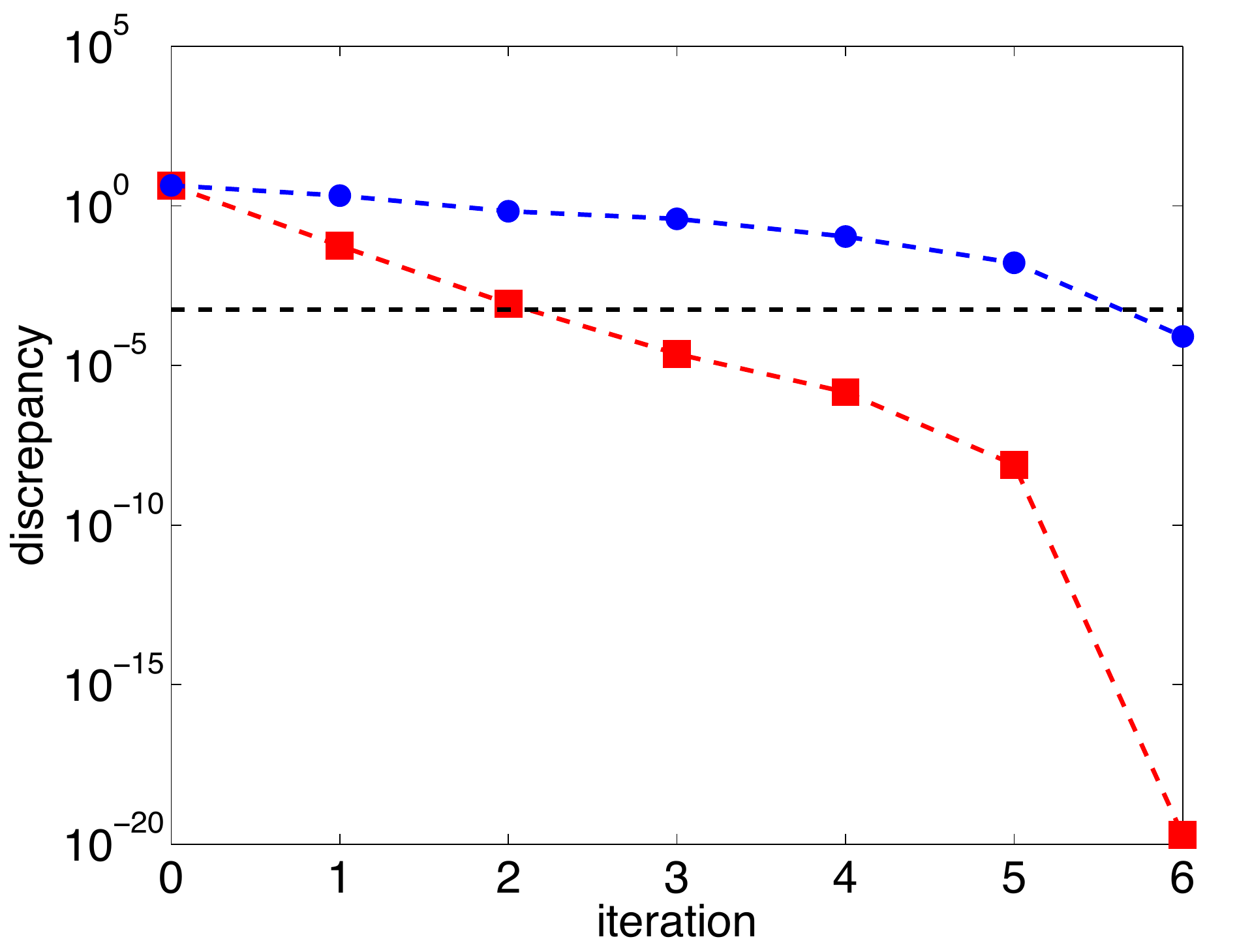}}
\caption{\label{fig:residual deconv} The residual of the iterative solutions as a function of the iteration number. The red curve with square markers corresponds to the plain CGLS, while the preconditioned residuals are marked by blue dots. The noise level used in the stopping rule is indicated by the horizontal line.}
\end{figure}

Since the convergence rate of the CGLS method depends on the condition number of the coefficient matrix, on the basis of the above observation that the width of the interval defined by the nonzero eigenvalues of $\mA^{\mT}\mA$ is smaller than the corresponding interval for the priorconditioned problem, we expect that the plain CGLS method satisfies the stopping criterion in fewer iterations than required by the PCGLS method. Figure~\ref{fig:residual deconv} confirms that this is indeed the case: the plain CGLS method stops after three iterations, while six iterations are required for the PCGLS method to terminate.

\subsection{Example 2: Computerized tomography with few radiographs}

In the second computed example, we consider the problem of estimating a two-dimensional density map in a bounded domain from noisy observation of integrals over a sparse set of lines across the domain. For the sake of definiteness, let $\Omega = [-1/2,1/2]\times[-1/2,1/2]$ represent the domain of interest, and let $\rho:\Omega\to\R$ denote a piecewise continuous non-negative density function. For simplicity, we assume that ${\rm supp}(\rho)\subset B(0,1/2)$, the disc centered at the origin and with radius $1/2$. Let $\ell(s,\theta)$ denote a line segment crossing the domain $\Omega$,
\[
 \ell(s,\theta) = \big\{ y = (y_1,y_2) \in \Omega\mid y_1 \cos\theta + y_2\sin \theta = s\big\},\quad |s|<1/2,\; 0\leq\theta<2\pi.
\]

The data consist of noisy observations of the integrals
\begin{equation}\label{beer}
 g(s,\theta) = \int_{\ell(s,\theta)} \rho(y) dS(y) = \int_0^S \rho(y(t)) dt,
\end{equation}
where $y=y(t)$ indicates the parametrization of the line segment with respect to the arc length $t$ such that $y(0),y(S)\in\partial\Omega$.  The integrals are supported over a discrete set of lines, parametrized as $\{\ell(s_1,\theta_1),\ldots,\ell(s_m,\theta_m)\}$. It follows from the Fourier slice theorem \cite{Natterer} that, in theory, the knowledge of the integrals over all lines crossing  $\Omega$ is enough to determine $\rho$ uniquely. In practice, if the set of lines is sparse, the problem becomes underdetermined and therefore ill-posed.

To discretize the problem, we divide the square $\Omega$ into $n\times n$ square pixels of equal size, denoted by $\Omega_j$, $1\leq j\leq n^2 = N$. The density $\rho$ is approximated by a piecewise constant function,  $\rho\big|_{\Omega_j} = x_j$, and the model (\ref{beer}) for noiseless data is approximated by
\[
g(s_k,\theta_k) \approx \frac 1N \sum_{j=1}^N |\ell(s_k,\theta_k)\cap\Omega_j| x_j,
\]
where $|\ell(s_k,\theta_k)\cap\Omega_j|$ denotes the length of the intersection of the line $\ell(s_k,\theta_k)$ and the $j$th pixel. Assuming additive noise $e\in\R^m$, the model leads to the linear observation equation $b = \mA x + e$ with obvious notations.

We consider a family of priors of Whittle-M\`{a}tern type, defined as follows: Let $\Delta_D$ denote the Dirichlet Laplacian over $\Omega$, that is,
\[
 {\mathscr D}(\Delta_D) = \big\{ u\in H^2(\Omega)\mid u\big|_{\partial \Omega}=0\big\}.
\]
We define the precision operator, the inverse of the covariance operator, through the formula
\[
 {\mathscr K} = -\Delta_D + \frac 1{\lambda^2} {\rm id},
\]
where the parameter $\lambda>0$ is a correlation length. This class of priors has been discussed extensively in recent articles \cite{Roininen1,Roininen2}, and it has been  shown that their discrete approximations are, in a certain sense, discretization invariant.

\begin{figure}
\centerline{\includegraphics[width=6cm]{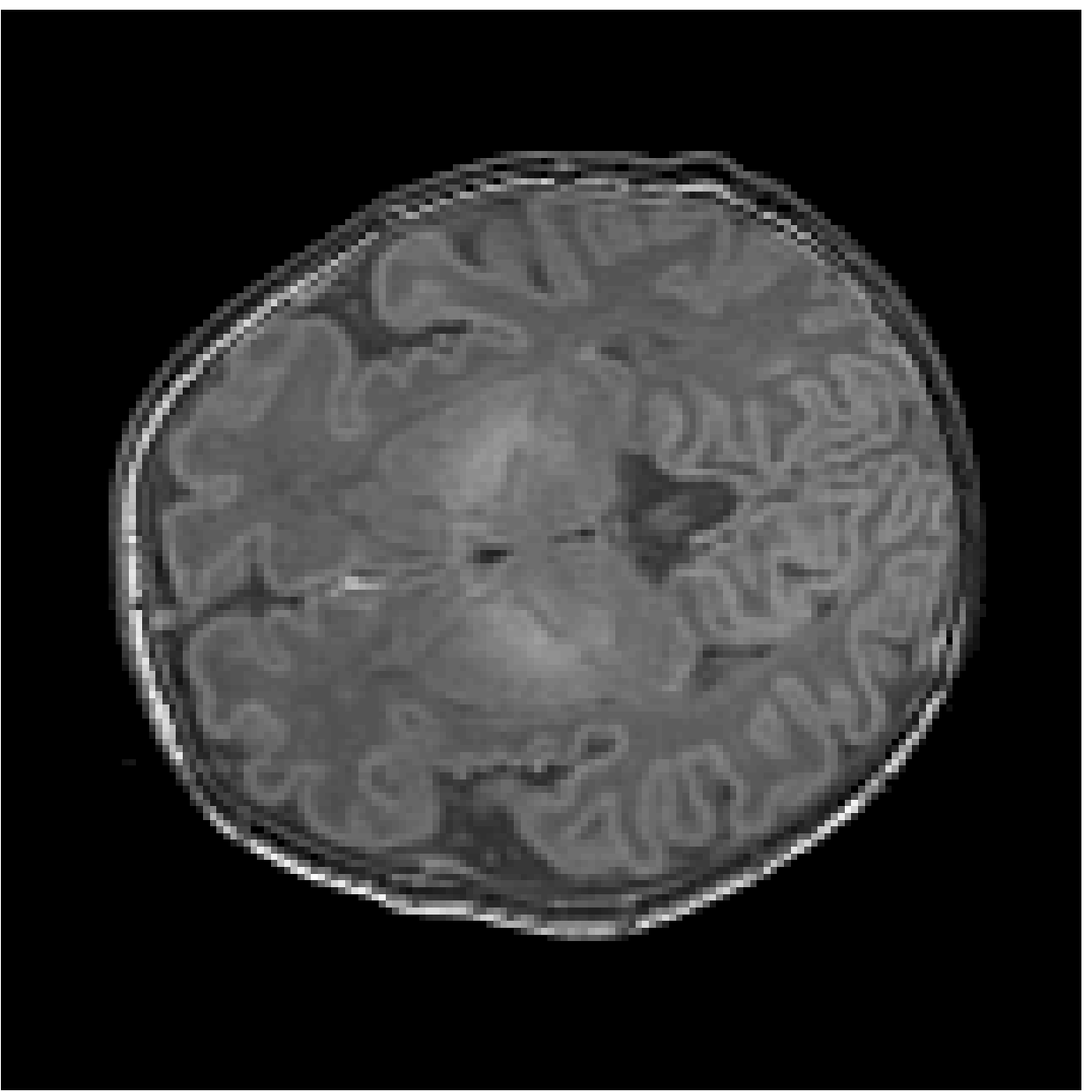}\includegraphics[width=7.8cm]{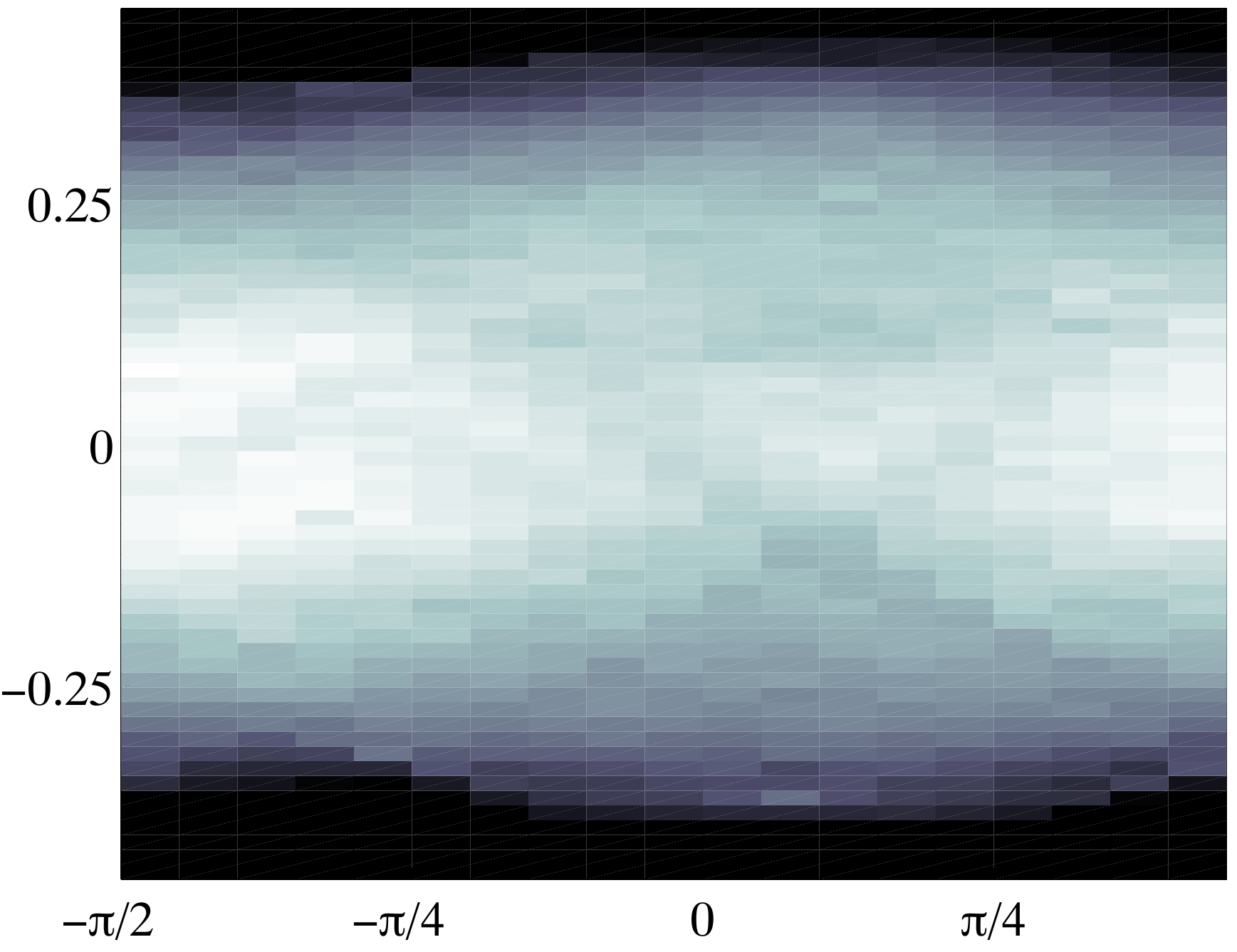}
}
\caption{\label{fig:true and data} The original image of size $160\times 160$ pixels (left) and the sinogram data (right), plotted as an image of size $60\times 20$. The horizontal axis corresponds to the 20 different illumination angles.}
\end{figure}

To discretize the prior, we use a finite difference approximation of the Laplacian, writing the precision matrix $\mK\in\R^{N\times N}$ as
\[
 \mK = -\mI_n\otimes \mD - \mD\otimes\mI_n + \frac 1{\lambda^2}\mI_N,
\]
where $\mD\in\R^{n\times n}$ is the three-point finite difference approximation of the one-dimensional  Laplacian with Dirichlet boundary conditions,
\[
 \mD = \frac 1{n^2} \left[\begin{array}{rrrr}    -2         &  1       &              & \\
                                                   1         & -2       & \ddots   & \\
                                                              &\ddots  &              &  1 \\
                                                               &            &    1        & -2 \end{array}\right],
\]
and $\otimes$ denotes the Kronecker product of matrices.
To generate the data, we consider the gray scale image of size $160\times 160$ shown in Figure~\ref{fig:true and data}, and illuminate it from  $n_\theta = 20$ illumination angles that are evenly distributed over the interval $[-\pi/2,\pi/2)$, $\theta_j = -\pi/2 + j \pi/n_\theta$, $0\leq j\leq n_\theta-1$.  For each illumination direction, $n_s=60$ parallel beams are chosen, corresponding to values $s_k = -1/2 + k/(n_s-1)$, $0\leq k  \leq n_s-1$, yielding to a severely undersampling of the sinogram data. The forward matrix obtained in this manner is of size $1\,200\times 25\,600$, hence the index of under determinacy is $n/m = 21.3$. The noiseless undersampled sinogram data is shown in Figure~\ref{fig:true and data}.

\begin{figure}
\centerline{\includegraphics[width=4cm]{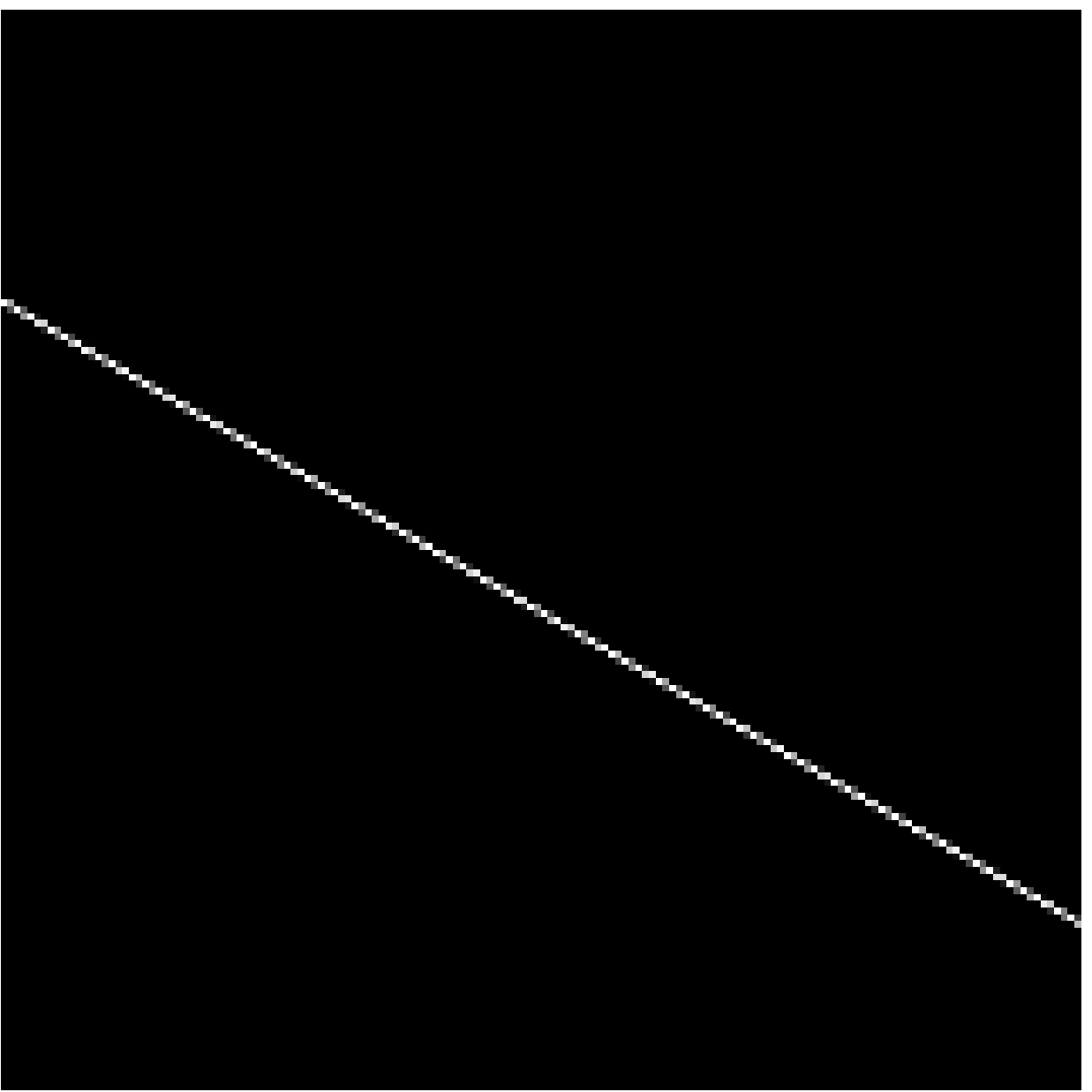}\includegraphics[width=4cm]{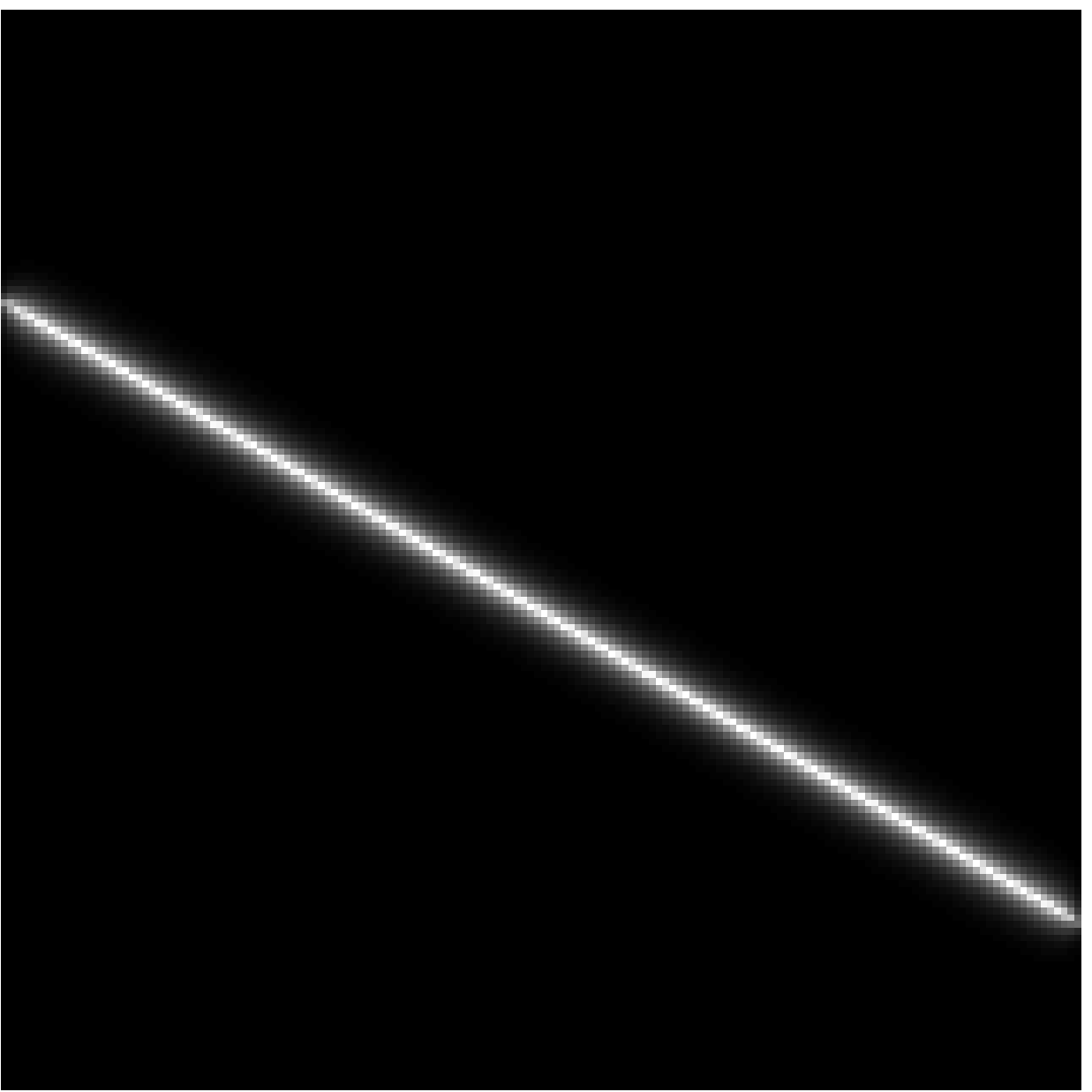}\includegraphics[width=4cm]{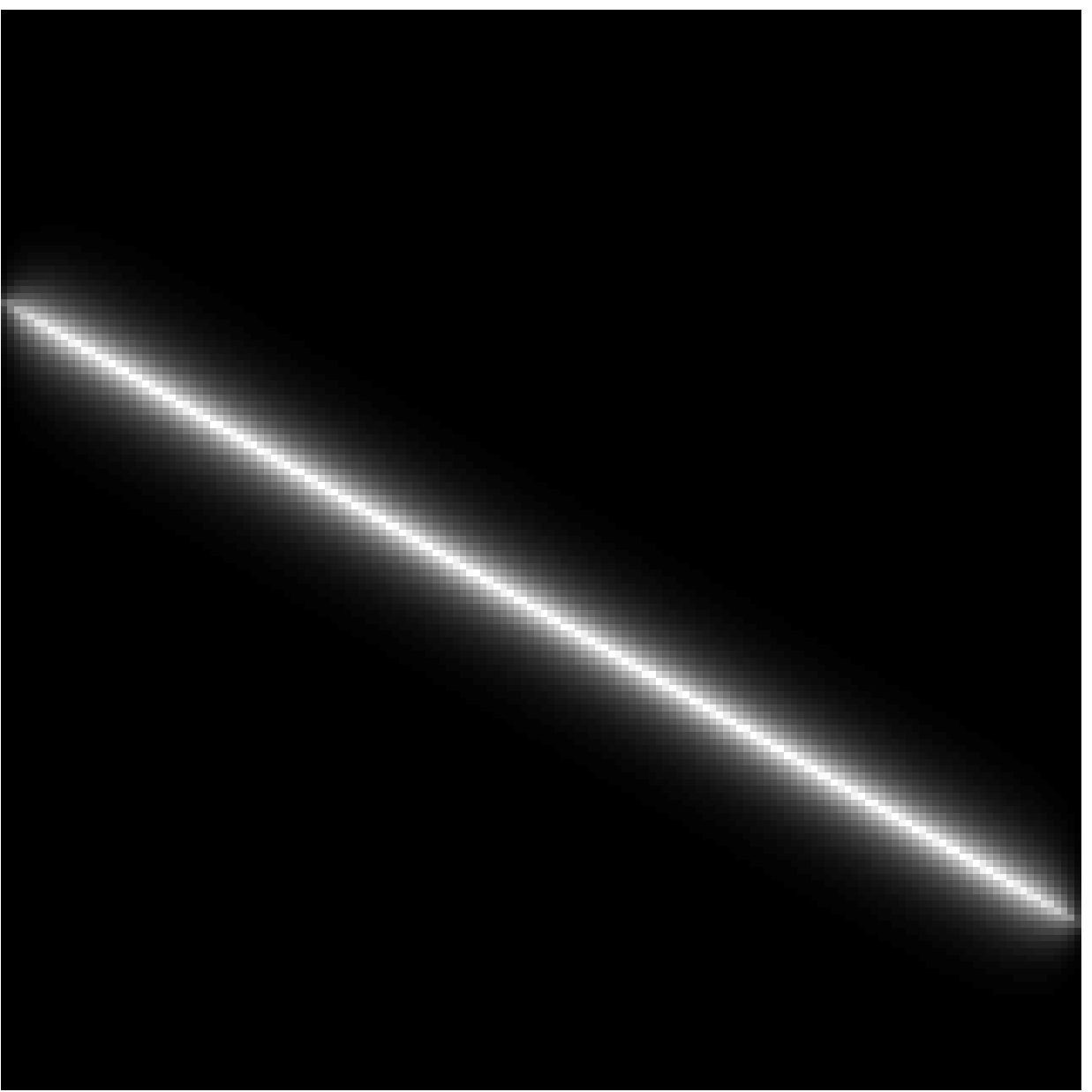}
}
\centerline{\includegraphics[width=4cm]{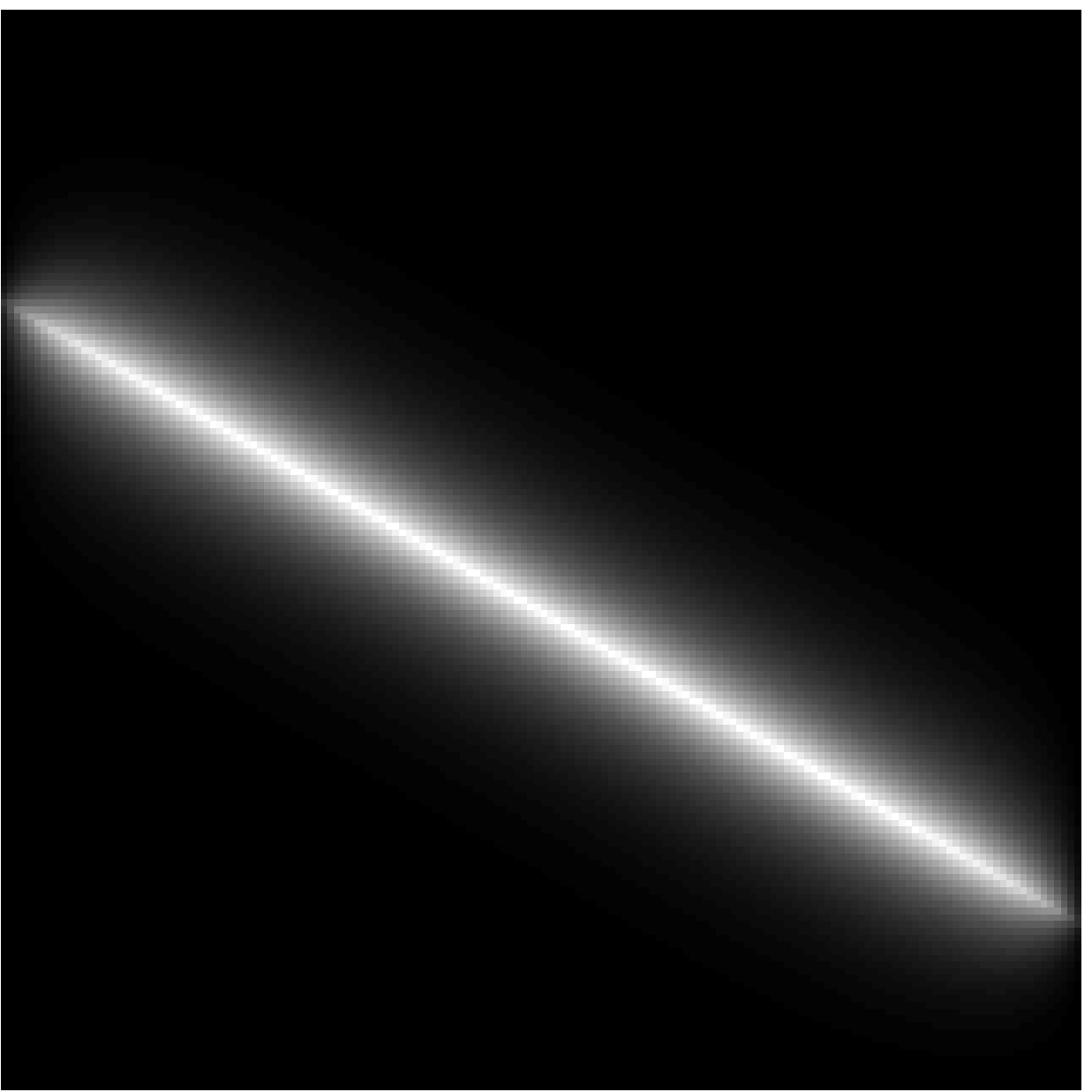}\includegraphics[width=4cm]{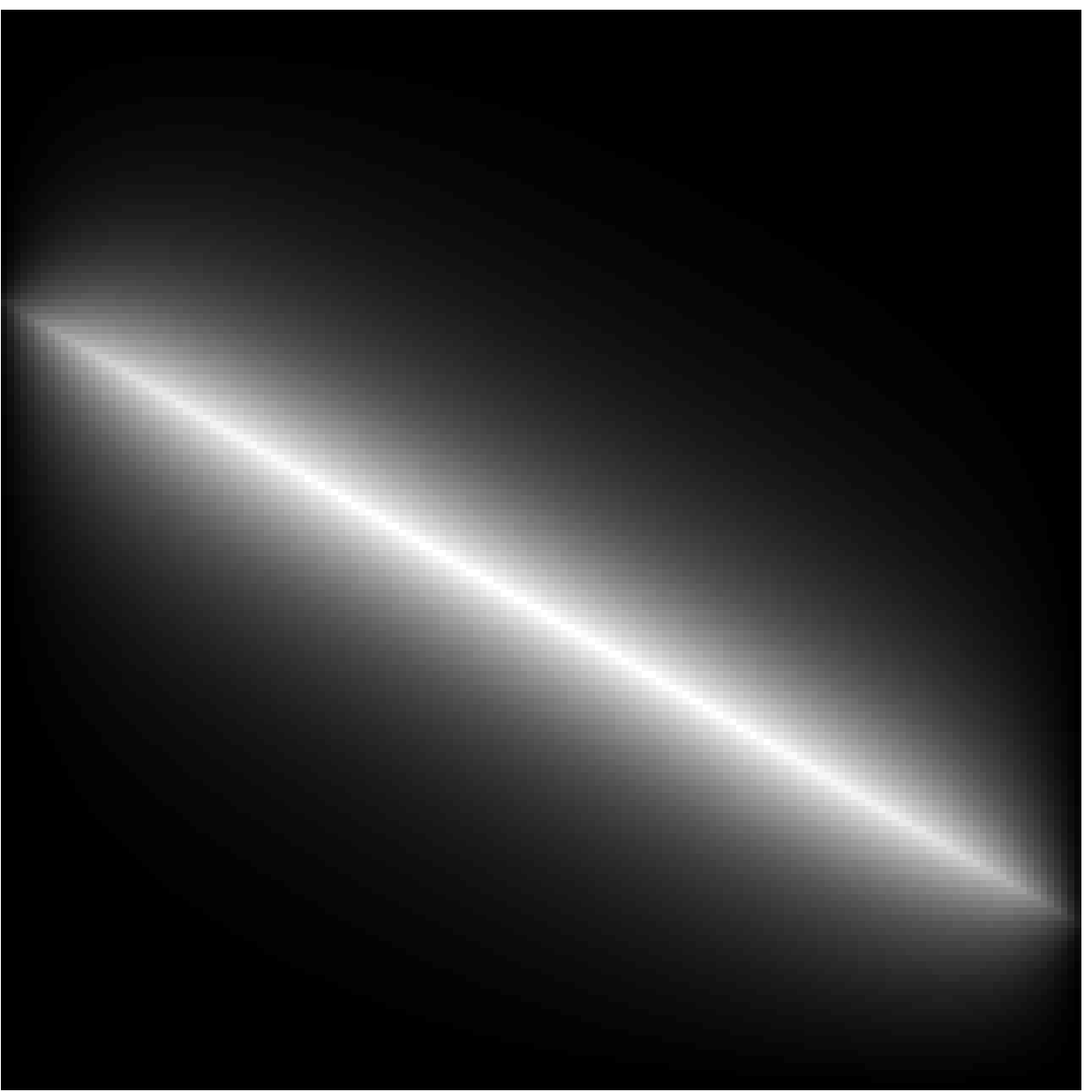}\includegraphics[width=4cm]{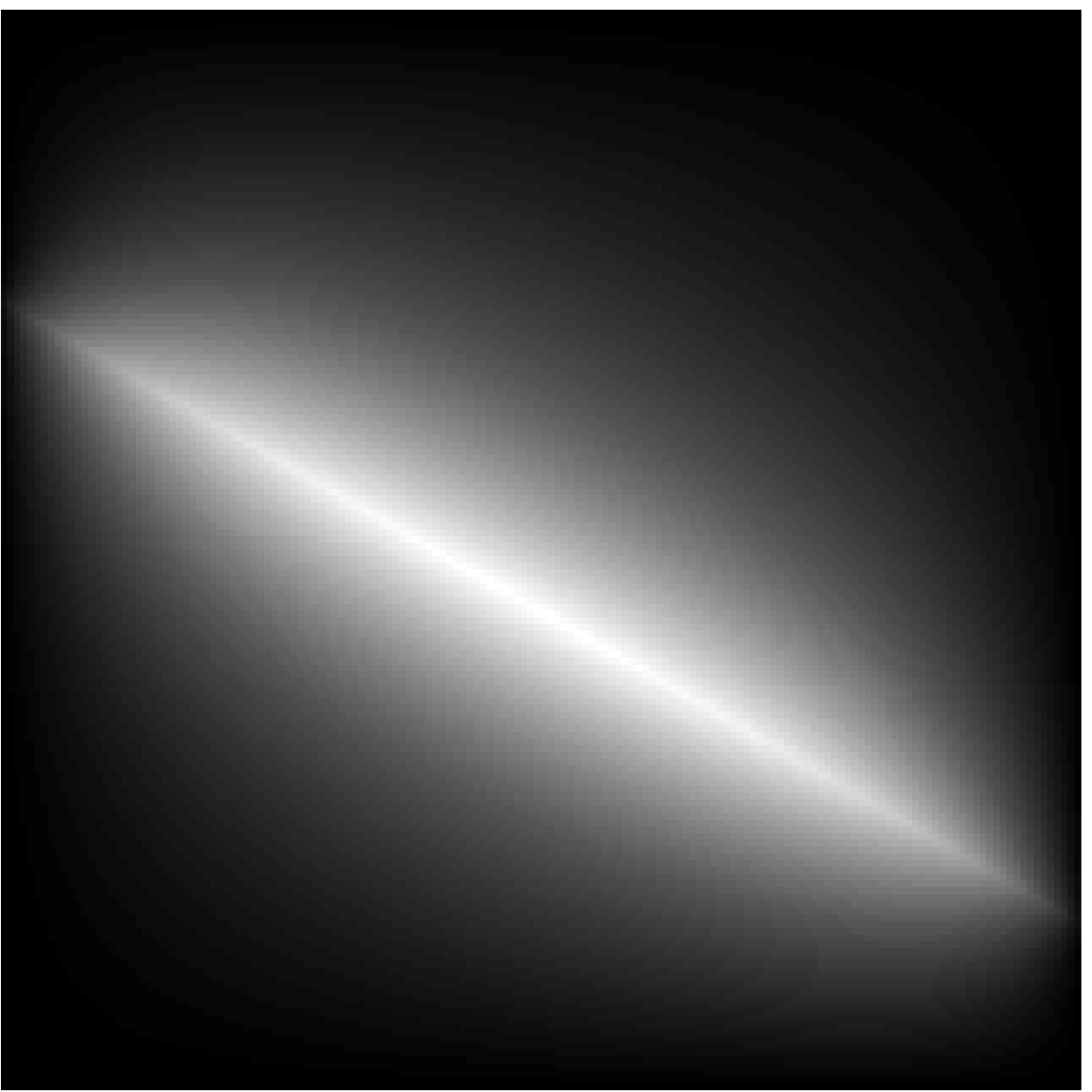}
}
\caption{\label{fig:basis vectors} A single column of the matrix ${\mA^\mT}$ representing a basis vector of the subspace ${\mathscr R}(\mA^\mT)$ represented as a $160\times 160$ image (top left), and the same vector after a multiplication of the covariance matrix $\mC = \mK^{-1}$, with different values of the correlation length parameter $\lambda$. Measured in units of pixels, the correlation length is chosen as $\lambda = 2,4$  pixels (top row), and $\lambda = 8,16$ and $32$ pixels (bottom row). }
\end{figure}

\begin{figure}
\centerline{\includegraphics[width=6cm]{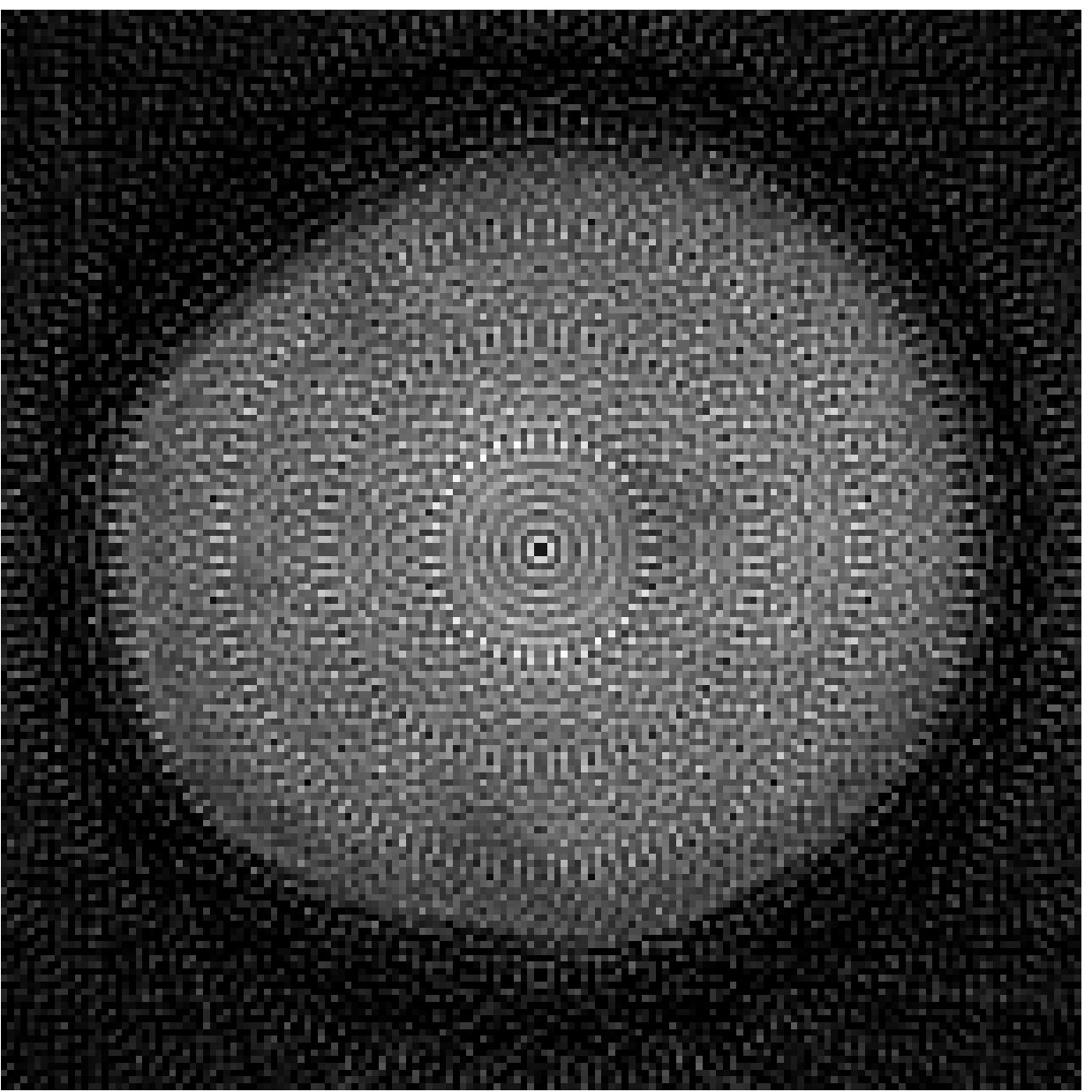}$\qquad$\includegraphics[width=6cm]{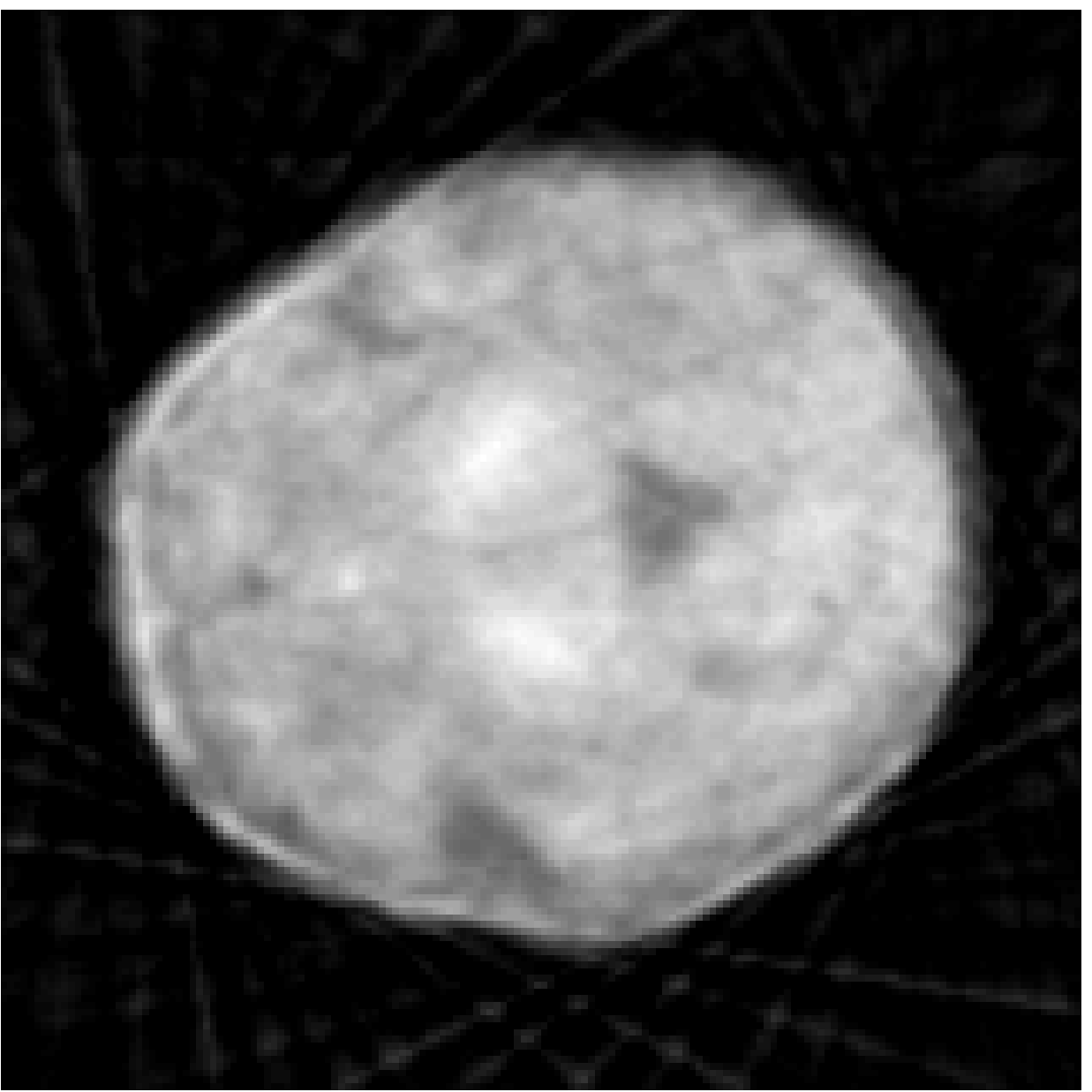}
}
\caption{\label{fig:reconstructions} The CGLS solutions with no priorconditioning (left) and with priorconditioning with correlation length $\lambda = 4$ pixels.}
\end{figure}

As in the previous example, we start by demonstrating the effect of priorconditioning on the basis vectors of ${\mathscr R}(\mA^\mT)$. In the present example, the column vectors of the matrix $\mA^\mT$ allow an interpretation as an image of size $160\times 160$. Therefore, in Figure~\ref{fig:basis vectors}, we show a column of the matrix, in image form, as well as the same vector multiplied by covariance matrices $\mC$ corresponding to different correlation lengths $\lambda$. As expected, the plain column vector represents a beam traversing the image area, while the multiplication with the covariance matrix creates blurring that increases with the correlation length. Observe that any image that is supported on pixels that have an empty intersection with all the beams defining the data are orthogonal to the columns of $\mA^\mT$, thus in the null space of $\mA$. Consequently, we expect that the plain CGLS algorithm produces an image that has a strong geometric artifact since the pixels with empty intersection with the beams remain dark. The application of the covariance operator on the basis vectors, on the other hand, increases the width of the beams, hence the preconditioned CGLS solution does not display this artifact, producing a smoother solution. This is indeed the case, as illustrated in Figure~\ref{fig:reconstructions}, where the two solutions produced by the CGLS method with and without prior conditioner are displayed side by side.

\begin{figure}
\centerline{\includegraphics[width=8cm]{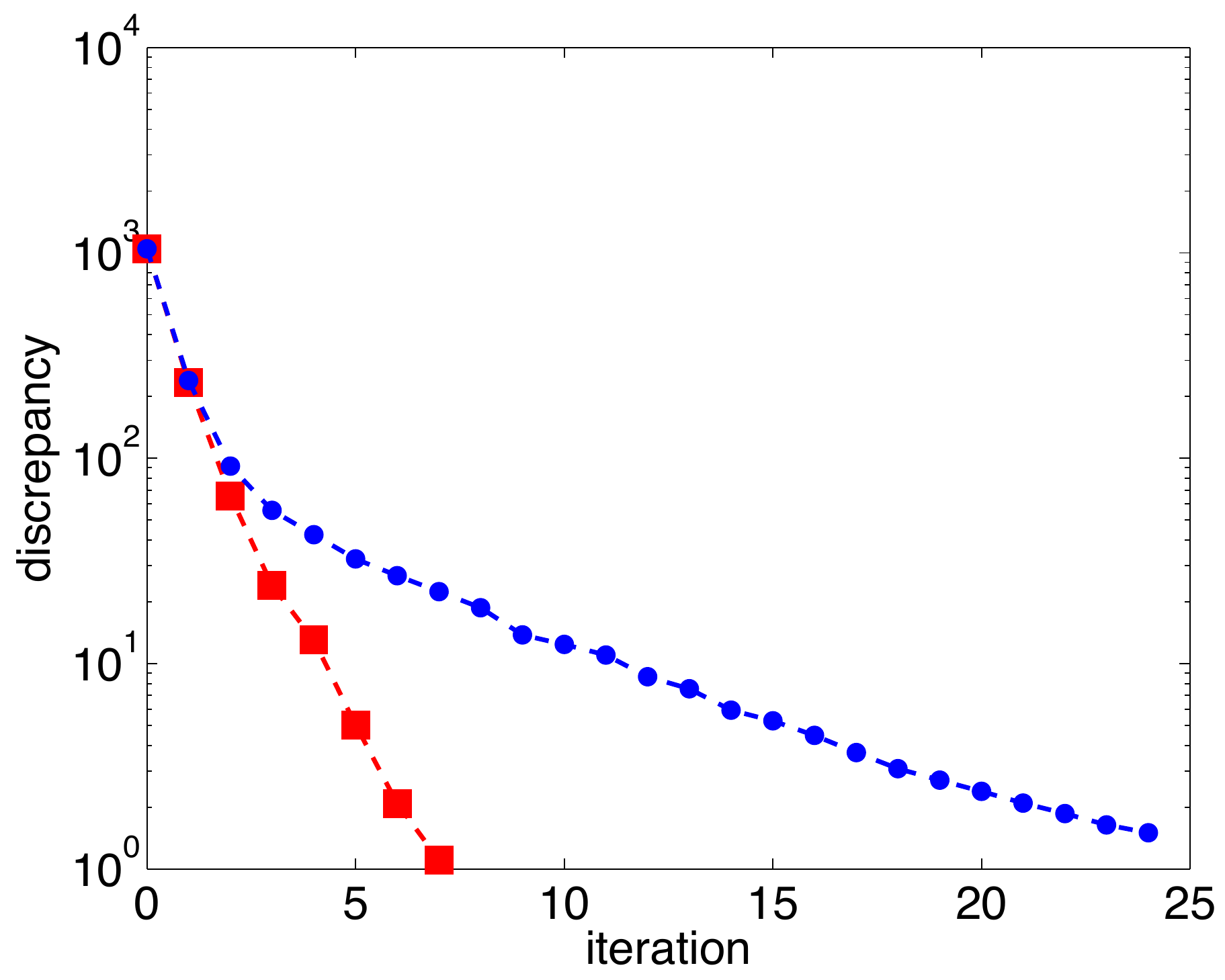}
\includegraphics[width=7.7cm]{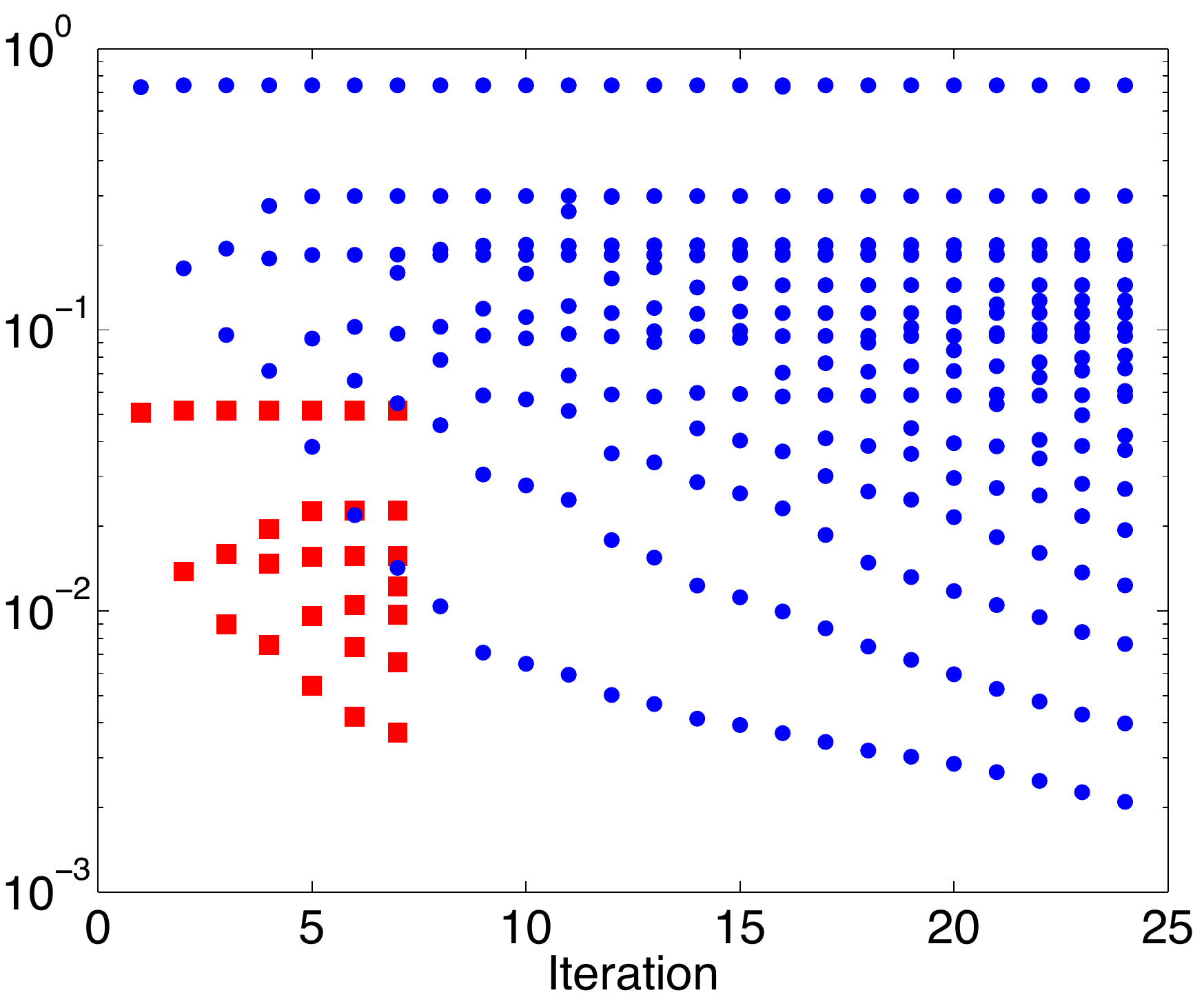}
}
\caption{\label{fig:discrepancy} Left: The discrepancy as a function of the iteration round, the plain CGLS corresponding to the red squares, and the prior conditioned CGLS to the blue dots. Right: The spectrum of the tridiagonal matrix $\mT_k$ as a function of iteration round. The red squares correspond to the plain CGLS, the blue dots to the priorconditioned CGLS. }
\end{figure}

As in the previous example, the enrichment  of the Krylov subspaces by vectors that lie in null space of the matrix $\mA$ slows down the convergence of method. In the left panel of Figure~\ref{fig:discrepancy}, we have plotted the discrepancy as a function of the iteration round both for plain CGLS and the prior conditioned version. Moreover, we plot the eigenvalues of the tridiagonal matrix $\mT_k$ approximating the eigenvalues of $\mA^\mT\mA$, as well as the corresponding approximations for the prior conditioned version of the algorithm. The eigenvalue approximations indicate that the eigenvalues of the prior conditioned matrix are significantly more spread out than without prior conditioning.  Numerical experiments indicate that by increasing the correlation length parameter $\lambda$ and thus widening the beams sounding the unknown density, the null space of $\mA$ has a stronger role in the reconstructions and more iterations are needed for convergence to discrepancy.

\begin{figure}
\centerline{\includegraphics[width=9.5cm]{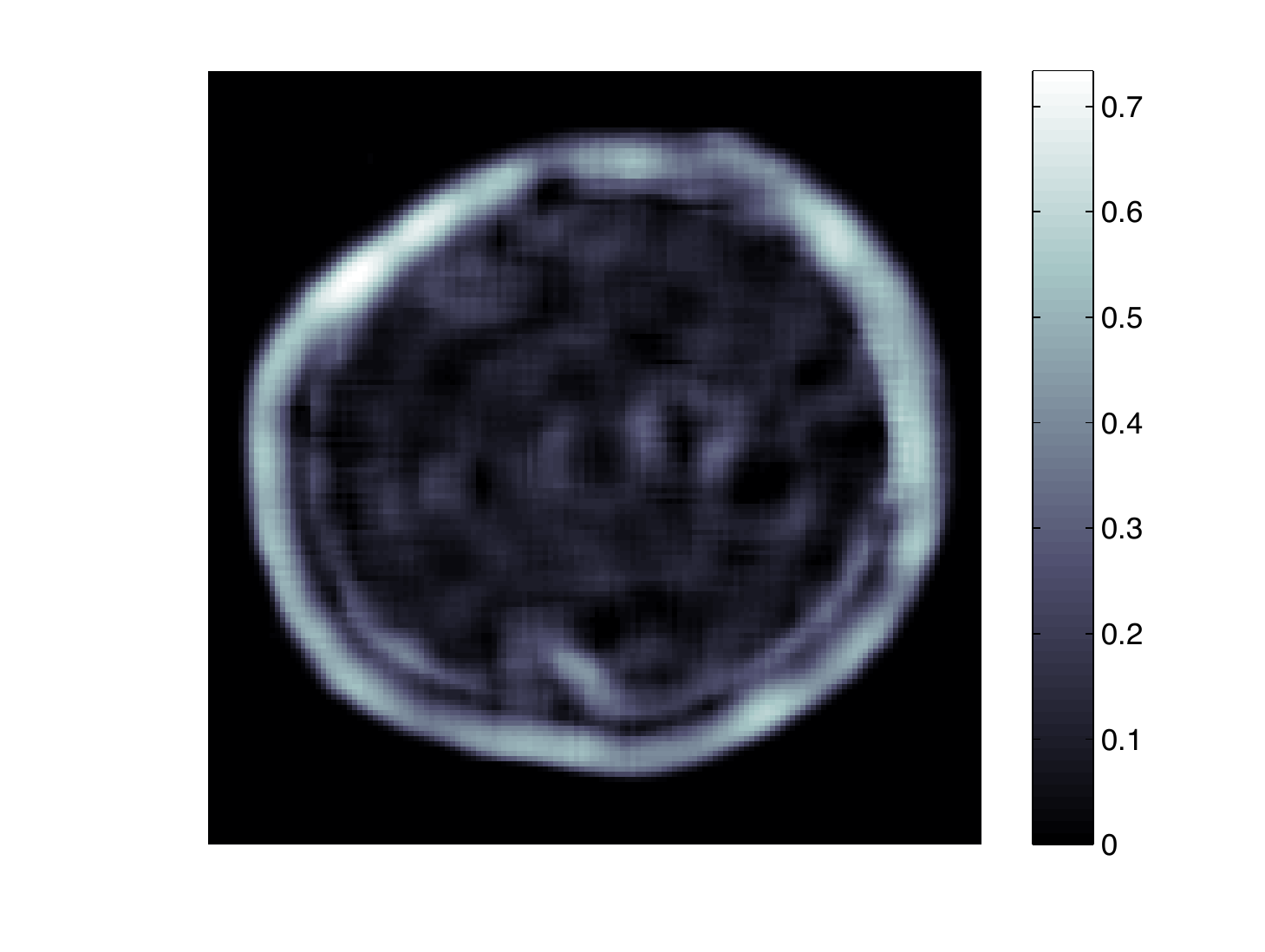} 
\includegraphics[width=9.5cm]{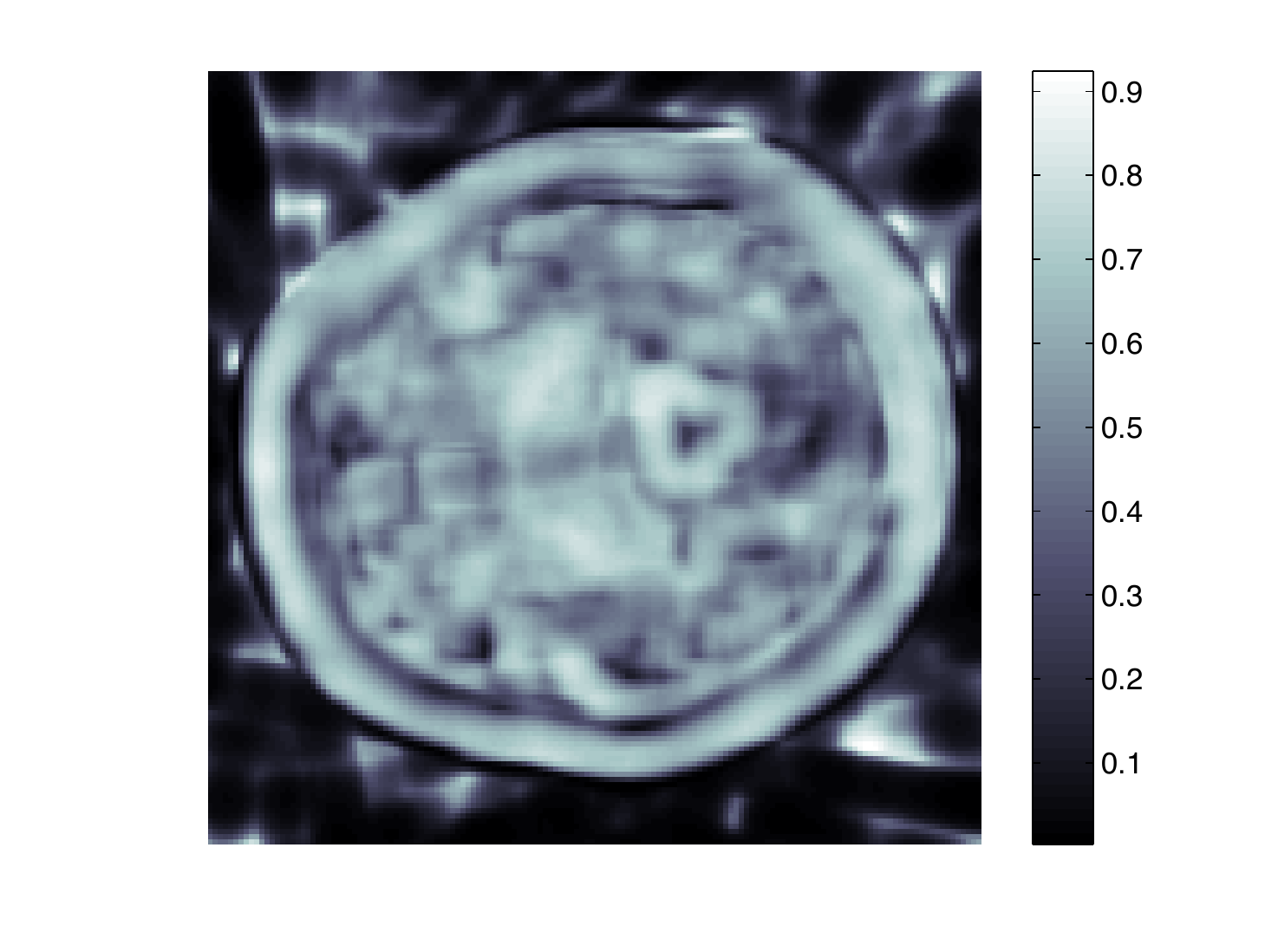}}
\caption{\label{fig:SSIM_maps} The similarity maps between the original image and the images obtained by CGLS with no priorconditioning (left) and with priorconditioning (right). Light colors denote the regions where the similarity is higher. }
\end{figure}

To evaluate the quality of the reconstructed images we use the Structural Similarity (SSIM) index
introduced in \cite{WBSS04}, which measures the structural differences between two images.
The SSIM index between the original image $I_{o}$ and the reconstructed
image $I_{rec}$ is defined as
$$
SSIM(I_{o},I_{rec}) = \frac{(2\,\mu_{o}\,\mu_{rec}+\gamma_1)(2\,\sigma_{o,rec}+\gamma_2)}
{(\mu_{o}^2+\mu_{rec}^2+\gamma_1)(\sigma_{o}^2+\sigma_{rec}^2+\gamma_2)}
$$
where $\mu_{o}$ ($\mu_{rec}$) and $\sigma_{o}$ ($\sigma_{rec}$) are the mean intensity and
the standard deviation of the original (reconstructed) image, respectively,
and  $\sigma_{o,rec}$ is the correlation between $I_{o}$ and $I_{rec}$.
We notice that the SSIM index is always $\le 1$, being equal to 1 when the images to be compared are identical.
\\
The two constants $\gamma_1$, $\gamma_2$ are added to avoid instability when the means $\mu_{o}$  and $\mu_{rec}$
or the standard deviation $\sigma_{o}$ and $\sigma_{rec}$ are small.
Typical values are $\gamma_1=0.01\,L$ and $\gamma_2=0.03\,L$, where $L$ is the dynamic
range of the images, i.e., the ratio between the largest and smallest values of the pixel intensities.
\\
In Figure~\ref{fig:SSIM_maps} the similarity maps between the original image (Figure~\ref{fig:true and data}, left)
and the reconstructed images (Figure~\ref{fig:reconstructions}) are shown. The maps are obtained by evaluating the local SIMM index
in a number of local Gaussian windows (see \cite{WBSS04} for details). The higher quality of the reconstruction
obtained by the solution with preconditioning is emphasized by values of the SSIM index
close to 1 in the region of interest. Compressing the similarity in a single indicator, the mean SSIM index is 0.120 when the solution is obtained
by the CGLS without priorconditioning, while its value is 0.409 when the priorconditioning is used.

\section*{Acknowledgements} This work was completed during the visit of DC and ES at University of Rome ``La Sapienza'' (Visiting Research/Professor Grant). The hospitality of the host university is kindly acknowledged. The work of ES was partly supported by NSF, Grant 1312424. This work was partially supported by grants from the Simons Foundation (\#305322  and \# 246665 to Daniela Calvetti).

\end{document}